\documentclass[10pt]{article}

\usepackage{dpr_fec,ifthen}
\usepackage{colortbl} 
\usepackage[graphicx]{realboxes} 
\usepackage{multirow} 
\usepackage{bigstrut} 
\usepackage{stackrel} 
\usepackage{cases}    
\hypersetup{colorlinks}
\allowdisplaybreaks

\usepackage{isetechreport}
\def\reportyear{18T}
\def\reportno{003}
\def\revisionno{1}
\def\originaldate{February 3, 2018}
\def\revisiondate{\today}
\coraltrue
\cvcrfalse

\begin{document}

\title{Regional Complexity Analysis of Algorithms for Nonconvex Smooth Optimization}

\author{Frank E.~Curtis\thanks{E-mail: \texttt{\href{mailto:frank.e.curtis@gmail.com}{frank.e.curtis@gmail.com}}}}
\affil{Department of Industrial and Systems Engineering, Lehigh University}
\author{Daniel P.~Robinson\thanks{E-mail: \texttt{\href{mailto:daniel.p.robinson@gmail.com}{daniel.p.robinson@gmail.com}}}}
\affil{Department of Applied Mathematics and Statistics, Johns Hopkins University}
\titlepage

\maketitle

\begin{abstract}
  A strategy is proposed for characterizing the worst-case performance of algorithms for solving nonconvex smooth optimization problems.  Contemporary analyses characterize worst-case performance by providing, under certain assumptions on an objective function, an upper bound on the number of iterations (or function or derivative evaluations) required until a $p$th-order stationarity condition is approximately satisfied.  This arguably leads to conservative characterizations based on anomalous objectives rather than on ones that are typically encountered in practice.  By contrast, the strategy proposed in this paper characterizes worst-case performance separately over regions comprising a search space.  These regions are defined generically based on properties of derivative values.  In this manner, one can analyze the worst-case performance of an algorithm independently from any particular class of objectives.  Then, once given a class of objectives, one can obtain an informative, fine-tuned complexity analysis merely by delineating the types of regions that comprise the search spaces for functions in the class.  Regions defined by first- and second-order derivatives are discussed in detail and example complexity analyses are provided for a few fundamental first- and second-order algorithms when employed to minimize convex and nonconvex objectives of interest.  It is also explained how the strategy can be generalized to regions defined by higher-order derivatives and for analyzing the behavior of higher-order algorithms.
\end{abstract}

\newcommand{\finf}{f_{\textnormal{inf}}}
\newcommand{\fref}{f_{\textnormal{ref}}}
\newcommand{\region}[1]{\Rcal_{#1}}
\newcommand{\subregion}[2]{\Rcal_{#1}^{#2}}
\newcommand{\RG}{\textnormal{\texttt{RG}}}
\newcommand{\RGA}{\textnormal{\texttt{RG-A}}}
\newcommand{\TRG}{\textnormal{\texttt{TR-G}}}
\newcommand{\TRH}{\textnormal{\texttt{TR-H}}}
\newcommand{\RN}{\textnormal{\texttt{RN}}}
\newcommand{\RNA}{\textnormal{\texttt{RN-A}}}
\newcommand{\Rp}{\textnormal{\texttt{Rp}}}
\newtheorem{step}{Step}
\numberwithin{definition}{section}
\numberwithin{corollary}{section}
\numberwithin{theorem}{section}
\numberwithin{lemma}{section}

\section{Introduction}\label{sec.introduction}

Users of optimization algorithms often choose to employ one algorithm instead of another based on its theoretical properties.  One such property of broad interest is \emph{worst-case} complexity, wherein one measures the resources that an algorithm will require, in the worst case, to solve (approximately) a given problem.  In the context of convex optimization \cite{Nest04}, such worst-case complexity has for many years been stated in terms of an upper bound on the number of iterations (or function or derivative evaluations\footnote{For the sake of brevity, we focus on worst-case complexity in terms of upper bounds on the number of \emph{iterations} required until a termination condition is satisfied, although in general one should also take \emph{function} and \emph{derivative evaluation} complexity into account.  These can be considered in the same manner as iteration complexity in our proposed strategy.}) required until either the distance between an iterate and an element of the set of minimizers, measured with a suitable norm, is less than a threshold $\epsilon_x \in (0,\infty)$, or the difference between an iterate's objective value and the optimal objective value is less than a threshold $\epsilon_f \in (0,\infty)$.

In the context of nonconvex optimization, a similar strategy has been adopted.  However, since one generally cannot guarantee that a method for solving nonconvex optimization problems will produce iterates that converge to a global minimizer, or at least have corresponding objective values that converge to the global minimum, the common approach has been to determine a worst-case upper bound on the number of iterations until a $p$th-order stationarity measure is satisfied with error below a threshold $\epsilon_p \in (0,\infty)$.  For example, in a body of literature that has been growing in recent years (see, e.g., \cite{BirgGardMartSantToin17,BirgMart17,CartGoulToin11b,CartGoulToin17,CurtRobiSama17,CurtRobiSama17b,Duss17,DussOrba17,GoulPorcToin12,NestPoly06,RoyeWrig17}), the main measure of interest has been the number of iterations required until an algorithm is guaranteed to produce an iterate at which the norm of the gradient of the objective function---a first-order stationarity measure---is below $\epsilon_1 \in (0,\infty)$.

Unfortunately, when it comes to minimizing broad classes of nonconvex objective functions satisfying loose assumptions---such as only Lipschitz continuity of some low-order derivatives---these types of worst-case complexity guarantees are forced to take into account exceptional objectives such that, when an algorithm is employed to minimize them, its behavior might be considered atypical.  For example, in \cite{CartGoulToin10,CartGoulToin11c}, Cartis, Gould, and Toint show that the worst-case guarantees for a few well-known methods are tight, but this is done with objective functions that one can argue are not representative of those encountered in regular practice.

One might attempt to overcome this resulting discrepancy between theory and practice in various ways.  Some argue that it would be ideal to be able to characterize \emph{average-case} behavior of an algorithm rather than worst-case, such as has been studied for the simplex method for solving linear optimization problems; see \cite{Borg82,Smal83,SpieTeng04}.  However, it seems difficult to set forth a useful, valid, and widely accepted definition of an average case when minimizing nonconvex objectives, even if one restricts attention to a small class of functions of interest.  Alternatively, one might consider analyzing the behavior of algorithms separately when they are employed to minimize functions in different classes.  However, this approach to worst-case performance guarantees limits itself to certain classes of objectives.

The purpose of this paper is to propose a strategy for characterizing the worst-case performance of algorithms for solving nonconvex smooth optimization problems.  In order to offer a characterization both $(i)$ within contexts seen in typical practice and $(ii)$ without limiting attention to specific problem classes, we propose that an algorithm's behavior can be characterized using a \emph{regional complexity analysis} (\emph{RC analysis}, for short) that involves the following two steps.

\benumerate\label{enum.steps_preliminary}
  \item Given an algorithm, one can analyze its performance by characterizing the behavior that it would exhibit within different regions in a search space (as defined in this paper).  This involves quantifying the decrease in the objective function that can be guaranteed when the algorithm finds itself at (or near) a point at which an objective's derivative values satisfy certain generic properties.
  \item After Step~1 is complete, one can combine results for an algorithm over combinations of regions in order to derive informative, fine-tuned analyses that characterize the worst-case performance of the algorithm when it is employed to minimize a function for which the search space is covered by the combination of regions.  For example, if one combines the results for an algorithm corresponding to \emph{region~1} and \emph{region~2}, then one can derive a worst-case complexity bound for the algorithm when it is employed to minimize functions for which the corresponding search space is covered completely by \emph{regions 1 and 2}.  For the same algorithm, this might lead to a different complexity than for, say, functions for which the search space is covered by points in \emph{regions 1, 2, and 3}.
\eenumerate

An interesting way to motivate our proposed strategy is to consider the seminal work of Nesterov and Polyak in \cite{NestPoly06}.  In this work, given a particular algorithm (namely, a cubicly regularized Newton method) and a class of objective functions (e.g., star-convex or gradient-dominated functions), the authors show that the algorithm progresses through different \emph{phases} as it converges to a solution set.  As revealed by the analysis, whether the algorithm is in a particular phase depends on the difference between the objective function value at a given iterate and the optimal objective value.  Our characterization strategy differs from the approach in~\cite{NestPoly06} in important ways, most significantly in the way that we \emph{decouple the analysis of the algorithm from consideration of a particular class of functions}.  Rather than start with a class of functions, we start with generically defined regions with which one can analyze the performance of an algorithm using the steps above.  In this manner, one does not consider a class of functions until the analysis of the algorithm over a set of regions has been completed.  A benefit of our approach is that it leads to a consistent standard for comparing methods across different function classes.  This is demonstrated in this paper as we simultaneously analyze a set of first- and second-order methods (rather than only one, as in~\cite{NestPoly06}).  Also, by considering regions defined by second- and higher-order derivatives, our strategy allows one to consider classes of nonconvex objectives beyond those considered in~\cite{NestPoly06}.

Our strategy can also be seen as a more comprehensive approach than ones that can be found piecemeal in other recent papers.  (Not to say that our work subsumes all ideas from these other papers; they also discuss other issues not considered here.)  For example, in \cite{KariNutiSchm16} (resp.~\cite{CarmHindDuchSidf17}), the authors show how gradient descent (resp.~accelerated gradient descent) exhibits a fast rate of convergence, even when minimizing a nonconvex function, \emph{if} it happens to take a path through the search space along which the function exhibits properties as if it were (strongly) convex.  In the case of \cite{CarmHindDuchSidf17}, if this behavior is not exhibited, then it is shown that an alternative type of step can be computed that would be beneficial to follow.  These articles show that the behavior of an algorithm can be better than that revealed by a contemporary worst-case analysis in a nonconvex setting, although neither paper sets forth a strategy for analyzing other types of algorithms, as we do.  Our strategy is also more comprehensive than approaches taken by authors who have studied the performance of algorithms in neighborhoods about \emph{strict saddle} points and related concepts; see, e.g., \cite{dauphin2014identifying,GeHuanJinYuan15,JinGeNetrKakaJord17,LeePanPilSimJorRec17,LeeSimcJordRech16,liu2017noisy,PateMokhRibe17}.  Such neighborhoods are a subset of the second-order-derivative-based regions that we define in this paper.

\subsection{Contributions}

Our contributions relate to our proposed RC analysis for characterizing the performance of algorithms for solving nonconvex smooth optimization problems.  Benefits of our strategy and our related contributions can be summarized as the following.

\bitemize
  \item Our proposed RC analysis of the performance of a given algorithm can be performed \emph{independently} from any particular class of functions.
  \item Given a class of functions, an RC analysis can offer a more fine-tuned worst-case performance analysis than the contemporary approach that only considers the number of iterations until $p$th-order stationarity is attained (approximately).  RC analysis offers a straightforward way to answer common questions, such as ``What performance can I expect if I apply an algorithm for nonconvex optimization to solve a problem that happens to be convex?''  Such questions are not answered well by contemporary worst-case results, which are too conservative.
  \item We demonstrate the use of RC analysis for analyzing first-, second-, and higher-order algorithms when employed to minimize functions in various classes of interest.  By tying the definitions for regions to properties of derivative values, RC analysis appropriately reveals performance guarantees that are representative of what can be expected in practice by derivative-based algorithms.
  \item RC analysis can be used to guide the design of \emph{new algorithms}.  For example, as demonstrated in this paper, an adaptive algorithm that computes different types of steps depending on properties of derivative values at a given iterate can achieve better RC analysis results than an algorithm that is not adaptive.  By contrast, using the contemporary approach to worst-case performance analysis, one often finds that certain static algorithms---such as gradient descent with a fixed stepsize or a cubicly regularized Newton method with a fixed regularization parameter---are optimal with respect to worst-case performance \cite{CartGoulToin10,CartGoulToin11c} despite the fact that adaptive algorithms often perform better in practice.
\eitemize

\subsection{Notation}

We use $\R{}$ to denote the set of real numbers (i.e., scalars), $\R{}_{\geq0}$ (resp.,~$\R{}_{>0}$) to denote the set of nonnegative (resp.,~positive) real numbers,~$\R{n}$ to denote the set of $n$-dimensional real vectors, and $\R{m \times n}$ to denote the set of $m$-by-$n$-dimensional real matrices.  The set of natural numbers is denoted as $\N{} := \{0,1,2,\dots\}$.  We write $\lambda(M)$ to denote the left-most (with respect to the real line) eigenvalue of a real symmetric matrix $M$.  Given $a \in \R{}$, we define $(a)_- := \max\{0,-a\}$, which is a nonnegative scalar that is strictly positive if and only if $a$ is strictly negative.  All norms are considered Euclidean; i.e., we let $\|\cdot\| := \|\cdot\|_2$.

If $\{a_k\}$ and $\{b_k\}$ are sequences of nonnegative scalars (i.e., elements of $\R{}_{\geq0}$), then we write $a_k = \Ocal(b_k)$ to indicate that there exists a positive constant $c \in \R{}_{>0}$ such that $a_k \leq cb_k$ for all $k \in \N{}$.  On the other hand, we write $a_k = \Omega(b_k)$ to indicate that there exists $c \in \R{}_{>0}$ such that $a_k \geq cb_k$ for all $k \in \N{}$.

Our problem of interest is to minimize $f(x)$ with respect to $x \in \R{n}$.  For simplicity, we assume that $f$ is real-valued and that one is interested in analyzing the behavior of a (monotone) descent algorithm, i.e., one for which, given an initial point $x_0 \in \R{n}$, the sequence $\{f(x_k)\}$ is monotonically nonincreasing over $\Lcal := \{x \in \R{n} : f(x) \leq f(x_0)\}$.  (Our strategies could also be extended to situations in which $f$ is extended-real-valued and for analyzing nonmonotone methods; see~\S\ref{sec.conclusion}.)  We append a natural number as a subscript for a quantity to denote its value during an iteration of an algorithm; e.g., henceforth, we let $f_k := f(x_k)$.

We make the following Assumption~\ref{ass.first} throughout the paper and add Assumption~\ref{ass.second} when analyzing second-order methods.

\bassumption\label{ass.first}
  The function $f : \R{n} \to \R{}$ is continuously differentiable and bounded below by $\finf := \inf_{x\in\R{n}} f(x) \in \R{}$.  In addition, over some open convex set $\Lcal^+$ containing $\Lcal$, the gradient function $g := \nabla f : \R{n} \to \R{n}$ is bounded in norm by $M_1 \in \R{}_{>0}$ and Lipschitz continuous with Lipschitz constant $L_1 \in \R{}_{>0}$; i.e.,
  \bequationNN
    \|g(x)\| \leq M_1\ \ \text{and}\ \ \|g(x) - g(\xbar)\| \leq L_1 \|x - \xbar\|\ \ \text{for all}\ \ (x,\xbar) \in \Lcal^+ \times \Lcal^+.
  \eequationNN
\eassumption

\bassumption\label{ass.second}
  Along with the conditions in Assumption~\ref{ass.first}, the function $f : \R{n} \to \R{}$ is twice continuously differentiable and, over the set $\Lcal^+$ defined in Assumption~$\ref{ass.first}$, the Hessian function $H := \nabla^2 f : \R{n} \to \R{n \times n}$ is Lipschitz continuous with Lipschitz constant $L_2 \in \R{}_{>0}$.  With Lipschitz continuity of $g$ from Assumption~$\ref{ass.first}$, the Hessian function is bounded in norm over $\Lcal^+$ by $M_2 \in \R{}_{>0}$, meaning that, overall,
  \bequationNN
    \|H(x)\| \leq M_2\ \ \text{and}\ \ \|H(x) - H(\xbar)\| \leq L_2 \|x - \xbar\|\ \ \text{for all}\ \ (x,\xbar) \in \Lcal^+ \times \Lcal^+.
  \eequationNN
\eassumption

As needed, specifically in \S\ref{sec.higher-order}, assumptions pertaining to higher-order continuous differentiability of~$f$ and Lipschitz continuity of higher-order derivatives of~$f$ will be introduced.  (Our strategies might also be extended to situations in which $f$ is nonsmooth and/or situations where the optimization problem involves implicit or explicit constraints.  We discuss such possibilities in~\S\ref{sec.conclusion}.)

\subsection{Algorithms}\label{sec.algorithms}

We will analyze the performance of a few algorithms throughout the paper.  It is important to note that this is done for demonstrative purposes only and that RC analysis is not limited to the types of algorithms considered here.  It can be employed to analyze the performance of other algorithms with different properties than those possessed by the algorithms that we discuss.  We comment on how such other analyses can be performed during our discussions.


\paragraph{Regularized gradient methods.}  We analyze two first-order methods, one static and one adaptive.  We refer to the static algorithm as the \emph{regularized gradient}~(\RG) method.  (It is often simply called gradient descent.  We call it \emph{regularized} for consistency in our naming scheme.)  At any iterate~$x_k$, this method produces the subsequent iterate as $x_{k+1} \gets x_k + s_k$, where, with $l_1 \in (L_1,\infty)$, one sets
\bequationNN
  s_k \gets \arg\min_{s\in\R{n}} f_k + g_k^Ts + \frac{l_1}{2} \|s\|^2 \implies x_{k+1} \gets x_k - \frac{1}{l_1} g_k.
\eequationNN
A similar, but adaptive first-order method, which we refer to as the \emph{adaptive regularized gradient}~(\RGA) method, computes a trial step at $x_k$ as $s_k \gets -g_k/\nu_k$ for some $\nu_k \in \R{}_{>0}$.  If this step yields a reduction in $f$ that is proportional to the reduction that it yields in the model $f_k + g_k^Ts + (\nu_k/2)\|s\|^2$, i.e.,
\bequation\label{eq.suff_dec_1}
  f_k - f(x_k+s_k) \geq \eta \(-g_k^Ts_k - \frac{\nu_k}{2} \|s_k\|^2\) = \frac{\eta}{2\nu_k} \|g_k\|^2
\eequation
for some $\eta \in (0,1)$, then the algorithm accepts the step by setting $x_{k+1} \gets x_k + s_k$; otherwise, it rejects it and $x_{k+1} \gets x_k$.  As for setting $\{\nu_k\}$, for $k = 0$ and any $k \geq 1$ such that $x_k \neq x_{k-1}$, the value $\nu_k$ is chosen from an interval $[\nu_{\min},\nu_{\max}] \subset \R{}_{>0}$; otherwise, if $s_k$ is rejected, then the method sets $\nu_{k+1} \gets \psi \nu_k$ for some $\psi \in (1,\infty)$.

\paragraph{Second-order trust region methods.}  Our next two algorithms are adaptive second-order trust region methods for which each trial step is computed as
\bequation\label{prob.tr}
  s_k \in \arg\min_{s\in\R{n}} f_k + g_k^Ts + \half s^TH_k s\ \ \text{subject to}\ \ \|s\| \leq \delta_k.
\eequation
The two methods that we consider merely differ in the manner in which $\{\delta_k\}$ is determined.  Both were studied in \cite[\S2.3--\S2.4]{CurtLubbRobi18}.  In the method we refer to as~\TRG, we let $\delta_k \equiv \|g_k\|/\nu_k$.  In the method we refer to as~\TRH, we let
\bequationNN
  \delta_k \equiv \frac{1}{\nu_k} \bcases \|g_k\| & \text{if $\|g_k\|^2 \geq (\lambda(H_k))_-^3$} \\ (\lambda(H_k))_- & \text{otherwise.} \ecases
\eequationNN
For \TRG{} and \TRH, $\{\nu_k\}$ is determined as in the \RGA{} method (for simplicity, using the same $\eta \in (0,1)$ and $\psi \in (1,\infty)$), except that in place of \eqref{eq.suff_dec_1} the methods observe
\bequation\label{eq.suff_dec_tr}
  f_k - f(x_k+s_k) \geq \eta\(-g_k^Ts_k - \half s_k^TH_ks_k\),
\eequation
which compares the reduction that the step offers in $f$ to the reduction that it offers in the second-order model $f_k + g_k^Ts + (1/2)s^TH_ks$.

\paragraph{Regularized Newton methods.}  We also consider two other second-order algorithms, but with different properties than the trust region methods stated above.  The first, a static second-order method that we refer to as the \emph{regularized Newton} (\RN) method, uses the update $x_{k+1} \gets x_k + s_k$, where, with $l_2 \in (L_2/2,\infty)$, it computes
\bequation\label{eq.RN}
  s_k \in \arg\min_{s\in\R{n}} f_k + g_k^Ts + \half s^TH_ks + \frac{l_2}{3} \|s\|^3.
\eequation
A similar, but adaptive method, which we refer to as the \emph{adaptive regularized Newton} (\RNA) method, computes trial steps as in \eqref{eq.RN}, but with $l_2$ replaced by $\nu_k$.  The sequence~$\{\nu_k\}$ is determined as in the \RGA, \TRG, and \TRH{} methods, except that for the \RNA{} method the employed sufficient decrease condition is
\bequationNN
  f_k - f(x_k + s_k) \geq \eta\(-g_k^Ts_k - \half s_k^TH_ks_k - \frac{\nu_k}{3} \|s_k\|^3\),
\eequationNN
which compares the reduction that the step offers in $f$ with the reduction that it offers in the regularized second-order model $f_k + g_k^Ts + (1/2)s^TH_ks + (\nu_k/3)\|s\|^3$.  (Again, we let \RNA{} use the same prescribed $\eta \in (0,1)$ and $\psi \in (1,\infty)$.)

We analyze the performance of other methods along with our discussion of higher-order RC analysis in \S\ref{sec.higher-order}.  We leave our description of those methods and the notation needed to state them for that section.

The algorithms described above as well as various other similar methods have appeared in the literature; see, e.g., \cite{BirgGardMartSantToin17,CartGoulToin11a,CartGoulToin11b,ConnGoulToin00,CurtLubbRobi18,FanYuan01,grapiglia2015convergence,grapiglia2016nonlinear,gratton2017decoupled,NestPoly06,NoceWrig06,Toin13}.  Therefore, for convenience, we draw from the literature when certain properties of these methods are needed.  That said, let us emphasize that critical aspects of our analyses of these methods are new, offering new perspectives on their behavior.

\subsection{Organization}

In \S\ref{sec.first-order}, we define regions based on first-order derivatives for our RC analysis framework, then analyze the behavior of the methods from~\S\ref{sec.algorithms} when an iterate lies in these regions.  We continue in \S\ref{sec.second-order} to define regions based on second-order derivatives, then analyze the performance of these algorithms when an iterate lies in these regions.  In~\S\ref{sec.summary}, we summarize our RC analysis results for these first- and second-order algorithms and provide complete perspectives on their behavior when minimizing functions in a few classes of interest.  Further discussion about, and possible variations on, the results in \S\ref{sec.summary} are presented in~\S\ref{sec.discussion}.  In~\S\ref{sec.higher-order}, we show how our framework can be generalized to regions defined according to higher-order derivatives and to analyze methods that employ higher-order derivative information.  Concluding remarks and ideas for extending RC analysis to other settings are provided in \S\ref{sec.conclusion}.

\section{First-Order Regions: Points with Gradient Domination}\label{sec.first-order}

We start to introduce our notion of regions with the following definition.  For this definition, recall that the first-order necessary condition for stationarity with respect to a continuously differentiable function $f$ is that $g(x)=0$.

\bdefinition[Region $\region{1}$]\label{def.1}
  \textit{
  For an objective $f : \R{n} \to \R{}$, scalar $\kappa \in (0,L_1]$, and reference objective value $\fref \in [\finf,\infty)$, let 
  \bequation\label{eq.gradient-dominated}
    \region{1} := \{x\in \Lcal: \|g(x)\|^\tau \geq \kappa(f(x) - \fref) \geq 0\ \text{for some}\  \tau \in [1,2] \}.
  \eequation
  Further, let $\subregion{1}{2}$ be the subset of $\region{1}$ such that the inequality in \eqref{eq.gradient-dominated} holds with $\tau = 2$ and let $\subregion{1}{1} := \region{1}\setminus\subregion{1}{2}$ so that $\region{1} = \subregion{1}{1} \cup \subregion{1}{2}$ with $\subregion{1}{1} \cap \subregion{1}{2} = \emptyset$.
  }
\edefinition

For flexibility in this definition, we have introduced $\fref \in [\finf,\infty)$.  Generally speaking, when analyzing the performance of a single algorithm, one can imagine~$\fref$ as a placeholder for the limiting value $\lim_{k\to\infty} f_k$, where the possibility of this value being strictly larger than $\finf$ might be inevitable due to nonconvexity of $f$.  On the other hand, if one can ensure---for a particular class of functions that will ultimately be considered---that the algorithm of interest will converge to global optimality, then one can consider the reference value to be $\fref = \finf$.  We discuss the role played by this value, and issues related to it, further in \S\ref{sec.discussion}.

Nesterov and Polyak \cite{NestPoly06} discuss a notion similar to that in Definition~\ref{def.1}; in particular, they refer to a function as \emph{gradient-dominated of degree $\tau$} if, for any $x \in \Lcal$, the inequality in \eqref{eq.gradient-dominated} holds for $\fref = \finf$ and some fixed $\tau \in [1,2]$.\footnote{Some authors take the term gradient-dominated to mean gradient-dominated of degree 2.  We do not take this meaning since, as seen in \cite{NestPoly06} and in this paper, functions that are only gradient-dominated of degree 1 offer different and interesting results.}  This range for~$\tau$ can be justified in various ways.  For one thing, $\tau \in (0,1)$ disproportionately weighs the norm of the gradient (as a measure of first-order stationarity) at points where it is small in norm.  On the other hand, allowing $\tau \in (2,\infty)$ would cause certain nice functions (such as strongly convex quadratics) not to have $\region{1} = \Lcal$, which would be undesirable.  We discuss in \S\ref{sec.summary} that certain well-known classes of functions---some convex and some nonconvex---have the property that $\region{1} = \Lcal$.  For example, this property holds for convex functions when $\Lcal$ is compact.

For an RC analysis pertaining to $\region{1}$, one is not restricting attention only to gradient-dominated functions.  Rather, by analyzing the behavior of algorithms with respect to $\region{1}$, one obtains results that are relevant for gradient-dominated functions \emph{as well as} for any nonconvex function for which points in a search space satisfy the inequality in \eqref{eq.gradient-dominated}, whether or not this includes the entire search space.  For example, for the function illustrated in Figure~\ref{fig.region1}, the region $\region{1}$ covers most of the search space, but not quite all of it.  This means that an RC analysis over~$\region{1}$ for a given algorithm will capture the worst-case performance of the algorithm over most of the domain, though it would not provide guarantees on the number of iterations it might spend in $\Lcal\setminus\region{1}$.  (For this, an analysis over a region defined according to higher-order derivatives might fill in the gap; see \S\ref{sec.second-order} and \S\ref{sec.higher-order}.)

\bfigure[ht]
  \centering
  \begin{minipage}[c]{0.35\textwidth}
    \includegraphics[width=\textwidth,clip=true,trim=45 25 50 25]{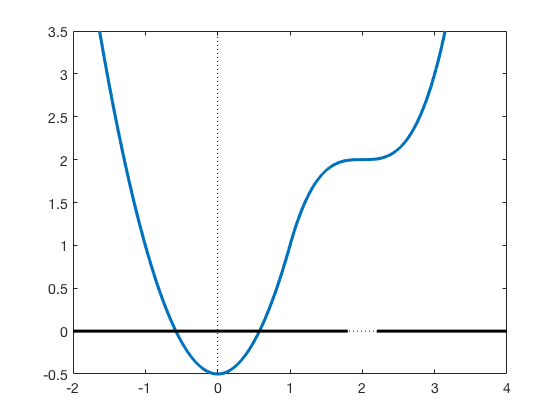}
  \end{minipage}\hfill
  \begin{minipage}[c]{0.63\textwidth}
    \caption[]{Plot of the continuously differentiable function
    \begin{minipage}{\linewidth}\bequationNN
      f(x) = \bcases \frac32 x^2 - \half & \text{if $x \leq 1$} \\ (x - 2)^3 + 2 & \text{otherwise.} \ecases
    \eequationNN \end{minipage}
    The bolded segments in the domain indicate $\region{1}$ for $\kappa=0.05$ and $\fref = \finf = -1/2$.  No matter the value for $\kappa \in \R{}_{>0}$, the region never includes an interval about $x = 2$, though it includes the rest of the domain.}
    \label{fig.region1}
  \end{minipage}
\efigure

Given this definition of $\region{1}$, one can provide insight into the performance of an algorithm merely by tying the reduction obtained with an accepted step to some gradient-related measure.  We formalize this with the following instruction, which should be understood as part of the first step introduced on page~\pageref{enum.steps_preliminary}.
\smallskip

\setcounter{step}{0}
\noindent
\fcolorbox{black}{gray!20}{\begin{minipage}{0.98\textwidth}
\begin{step}[Region $\region{1}$]\label{step.1_1}
  Attempt to prove that for any accepted step~$s_k$ the decrease in the objective function from $x_k$ to $x_{k+1} = x_k + s_k$ satisfies
  \bequation\label{eq.region1_reduction}
    f_k - f_{k+1} = \Omega(\|g(x)\|^r)\ \ \text{for some}\ x\in\{x_k,x_{k+1}\} \ \text{with} \  x\in\region{1} \ \text{and} \ r > 0.
  \eequation
  If such $(x,r)$ exists, then one can combine  \eqref{eq.gradient-dominated} and \eqref{eq.region1_reduction} to prove a reduction in the objective gap to~$\fref$, i.e., an upper bound for $f_{k+1} - \fref$ as a function of $f_k - \fref$.
  \end{step}
\end{minipage}}
\smallskip

It is implicit in~\eqref{eq.region1_reduction} that one considers the performance of an algorithm over~$\region{1}$ only when $\{x_k,x_{k+1}\} \cap \region{1} \neq \emptyset$.  This is reasonable since this is precisely when the size of the gradient at $x_k$ and/or $x_{k+1}$ gives information about the size of a potential reduction in the objective through the inequality~\eqref{eq.gradient-dominated} that defines $\region{1}$. 

In the remainder of this section, we provide two examples of following this instruction, which we refer to as Step~\ref{step.1_1}--$\region{1}$.  These will allow us to state results for the algorithms from~\S\ref{sec.algorithms}.  For our first theorem, we state a result pertaining to algorithms that, with an accepted step, yield a reduction in the objective that is proportional to the squared norm of the gradient at the current iterate.  This will allow us to characterize the behavior of \RG{}, \RGA{}, \TRG{}, and \TRH{} over $\region{1}$.

\btheorem\label{th.RC_region1_alg1}
  Suppose Assumption~$\ref{ass.first}$ holds.  For any algorithm such that $x_k\in\region{1}$ implies that~\eqref{eq.region1_reduction} holds with $x=x_k$ and $r = 2$ in that
  \bequation\label{eq.RC_region1_alg1}
    f_k - f_{k+1} \geq \frac{1}{\zeta} \|g_k\|^2\ \ \text{for some}\ \ \zeta \in [L_1,\infty),
  \eequation
  the following statements hold true.
  \benumerate
    \item[$($a$)$] If $x_k \in \subregion{1}{2}$, then $\{f_k - \fref\}$ decreases as in a linear rate; specifically,
    \bequation\label{eq.RC_region1_linear}
      f_{k+1} - \fref \leq \(1 - \frac{\kappa}{\zeta}\) (f_k - \fref)\ \ \text{where}\ \ \frac{\kappa}{\zeta} \in (0,1].
    \eequation
    \item[$($b$)$] If $x_k \in \subregion{1}{1}$, then it must be true that $\kappa(f_k - \fref) < 1$ and it follows that $\{f_k - \fref\}$ decreases as in a sublinear rate; specifically,
    \bequation\label{eq.RC_region1_sublinear}
      f_{k+1} - \fref \leq \(1 - \frac{\kappa^2}{\zeta}(f_k - \fref)\) (f_k - \fref).
    \eequation
  \eenumerate
  Similarly, for any algorithm such that having $x_k \in \region{1}$ implies that
  \bequation\label{eq.RC_region1_alg1_mstep}
    f_k - f_{k+m} \geq \frac{1}{\zeta} \|g_k\|^2\ \ \text{for some}\ \ \zeta \in [L_1,\infty)\ \ \text{and}\ \ m \in \N{}
  \eequation
  with $m$ independent of $k$, then $($a$)$ and $($b$)$ hold with $f_{k+1}$ replaced by $f_{k+m}$.\footnote{In this case, the decrease in the objective would be indicative of an \emph{$m$-step} linear (for part~(a)) or \emph{$m$-step sublinear} (for part~(b)) rate of convergence.  We do not explicitly refer to such a multi-step aspect of a convergence rate since it is always clear from the context.}
\etheorem
\bproof
  If $x_k \in \subregion{1}{2}$, then, with \eqref{eq.RC_region1_alg1}, it follows that
  \bequationNN
    f_k - f_{k+1} \geq \frac{1}{\zeta} \|g_k\|^2 \geq \frac{\kappa}{\zeta} (f_k - \fref).
  \eequationNN
  Adding and subtracting $\fref$ on the left-hand side and rearranging gives \eqref{eq.RC_region1_linear}.
  
  If $x_k \in \subregion{1}{1}$, which is to say that $\|g_k\| \geq \kappa (f_k - \fref)$ while $\|g_k\|^2 < \kappa (f_k - \fref)$, then it must be true that $\kappa(f_k - \fref) < 1$.  In this case, from \eqref{eq.RC_region1_alg1},
  \bequationNN
    f_k - f_{k+1} \geq \frac{1}{\zeta} \|g_k\|^2 \geq \frac{1}{\omega} (f_k - \fref)^2\ \ \text{where}\ \ \omega := \frac{\zeta}{\kappa^2}.
  \eequationNN
  Adding and subtracting $\fref$ on the left-hand side, one finds by defining the value $a_k := (f_k - \fref)/\omega = \kappa^2(f_k - \fref)/\zeta \in [0,1)$ for all $k \in \N{}$ that
  \bequationNN
    \underbrace{\frac{f_k - \fref}{\omega}}_{a_k} - \underbrace{\frac{f_{k+1} - \fref}{\omega}}_{a_{k+1}} \geq \underbrace{\frac{(f_k - \fref)^2}{\omega^2}}_{a_k^2}.
  \eequationNN
  One finds from this inequality that $a_{k+1} \leq (1 - a_k)a_k$, which gives \eqref{eq.RC_region1_sublinear}.
    
  If, with $x_k \in \region{1}$, an algorithm offers \eqref{eq.RC_region1_alg1_mstep}, then the desired conclusions hold using the same arguments above with \eqref{eq.RC_region1_alg1_mstep} in place of \eqref{eq.RC_region1_alg1}.
  \qed
\eproof

Not all algorithms offer inequality~\eqref{eq.RC_region1_alg1} (or even \eqref{eq.RC_region1_alg1_mstep}) while others offer an even stronger bound.  As our second example of following Step~\ref{step.1_1}--$\region{1}$, we prove the following theorem, which will allow us to characterize the behavior of our other second-order algorithms (i.e., \RN{} and \RNA{}) over $\region{1}$.  For the proof of this theorem, we use similar strategies as are used to prove \cite[Theorem~6 and Theorem~7]{NestPoly06}.  Interestingly, as for these results in \cite{NestPoly06}, one finds different behavior depending on whether $f_k - \fref$ is below a certain threshold.  For our purposes, we also need to consider a couple cases depending on properties of the iterates $x_k$ and $x_{k+1}$.\footnote{There arises an interesting scenario in this theorem for $x_{k+1} \in \subregion{1}{1}$ during which $\{f_k - \fref\}$ might initially decrease at a superlinear rate.  However, this should not be overstated.  After all, if this scenario even occurs, then the number of iterations in which it will occur will be limited if the iterates remain at or near points in $\region{1}$.}

\btheorem\label{th.RC_region1_alg2}
  Suppose Assumption~$\ref{ass.first}$ holds.  For any algorithm such that $x_{k+1} \in \region{1}$ implies that~\eqref{eq.region1_reduction} holds with $x = x_{k+1}$ and $r = 3/2$ in that
  \bequation\label{eq.RC_region1_alg2}
    f_k - f_{k+1} \geq \frac{1}{\zeta} \|g_{k+1}\|^{3/2}\ \ \text{for some}\ \ \zeta \in (0,\infty),
  \eequation
  the following statements hold true.
  \benumerate
    \item[$($a$)$] If $x_{k+1} \in \subregion{1}{2}$ and $f_k - \fref \geq \kappa^3/\zeta^4$, then $\{f_k - \fref\}$ has decreased as in a linear rate in the sense that the following inequality holds:
    \bequation\label{eq.cubic_region12_1}
      f_{k+1} - \fref \leq \(\frac{(f_0 - \fref)^{1/4}}{\frac{\kappa^{3/4}}{\zeta} + (f_0 - \fref)^{1/4}}\)(f_k - \fref).
    \eequation
    On the other hand, if $x_{k+1} \in \subregion{1}{2}$ and $f_k - \fref < \kappa^3/\zeta^4$, then the sequence has decreased as in a superlinear rate in the sense that
    \bequation\label{eq.cubic_region12_2}
      f_{k+1} - \fref \leq \(\frac{\zeta^4(f_k - \fref)}{\kappa^3}\)^{1/3} (f_k - \fref).
    \eequation
    \item[$($b$)$] If $x_{k+1} \in \subregion{1}{1}$, then it must be true that $\kappa(f_{k+1} - \fref) < 1$.  Thus, if $x_{k+1} \in \subregion{1}{1}$ and $f_k - \fref \geq \zeta^2/\kappa^3$, then $\{f_k - \fref\}$ has decreased as in a superlinear rate with
    \bequation\label{eq.cubic_region11_1}
      f_{k+1} - \fref \leq \(\frac{\zeta^2}{\kappa^3(f_k - \fref)}\)^{1/3}(f_k - \fref),
    \eequation
    whereas, if $x_{k+1} \in \subregion{1}{1}$ and $f_k - \fref < \zeta^2/\kappa^3$, then the sequence has decreased as in a sublinear rate in the sense that the following inequality holds:
    \bequation\label{eq.cubic_region11_2}
      f_{k+1} - \fref \leq \(\frac{1}{1 + \frac{\kappa^{3/2}}{\zeta}\(\frac{\sqrt{2}-1}{\sqrt{2}}\)\sqrt{f_k - \fref}}\)^2 (f_k - \fref).
    \eequation
  \eenumerate
  Similarly, for an algorithm such that having $x_{k+m} \in \region{1}$ implies that
  \bequation\label{eq.RC_region1_alg2_mstep}
    f_k - f_{k+m} \geq \frac{1}{\zeta} \|g_{k+m}\|^{3/2}\ \ \text{for some}\ \ \zeta \in (0,\infty)\ \ \text{and}\ \ m \in \N{}
  \eequation
  with $m$ independent of $k$, then $($a$)$ and $($b$)$ hold with the pair $(x_{k+1},f_{k+1})$ replaced by the pair $(x_{k+m},f_{k+m})$.
\etheorem
\bproof
  If $x_{k+1} \in \subregion{1}{2}$, then, with \eqref{eq.RC_region1_alg2}, it follows that
  \bequationNN
    f_k - f_{k+1} \geq \frac{1}{\zeta} \|g_{k+1}\|^{3/2} \geq \omega^{1/4} (f_{k+1} - \fref)^{3/4}\ \ \text{where}\ \ \omega := \frac{\kappa^3}{\zeta^4}.
  \eequationNN
  Adding and subtracting $\fref$ on the left-hand side, one finds by defining the values $a_k := (f_k - \fref)/\omega$ for all $k \in \N{}$ that
  \bequation\label{eq.3/4_a}
    \underbrace{\frac{f_k - \fref}{\omega}}_{a_k} - \underbrace{\frac{f_{k+1} - \fref}{\omega}}_{a_{k+1}} \geq \underbrace{\frac{(f_{k+1} - \fref)^{3/4}}{\omega^{3/4}}}_{a_{k+1}^{3/4}}.
  \eequation
  One finds from this inequality and monotonicity of $\{a_k\}$ that
  \bequationNN
    \frac{a_k}{a_{k+1}} \geq 1 + \frac{1}{a_{k+1}^{1/4}} \geq 1 + \frac{1}{a_0^{1/4}} \in (1,\infty),
  \eequationNN
  which gives \eqref{eq.cubic_region12_1}.  That said, if $a_k < 1$ (which is to say that $f_k - \fref < \omega = \kappa^3/\zeta^4$), then one finds from \eqref{eq.3/4_a} and $a_{k+1}\geq0$ that $a_{k+1} \leq a_k^{4/3}$, from which \eqref{eq.cubic_region12_2} follows.
    
  If $x_{k+1} \in \subregion{1}{1}$, which is to say that $\|g_{k+1}\| \geq \kappa (f_{k+1} - \fref)$ while $\|g_{k+1}\|^2 < \kappa (f_{k+1} - \fref)$, then it must be true that $\kappa(f_{k+1} - \fref) < 1$.  Hence, with \eqref{eq.RC_region1_alg2},
  \bequationNN
    f_k - f_{k+1} \geq \frac{1}{\zeta}\|g_{k+1}\|^{3/2} \geq \omega^{-1/2} (f_{k+1} - \fref)^{3/2}\ \ \text{where}\ \ \omega := \frac{\zeta^2}{\kappa^3}.
  \eequationNN
  Adding and subtracting $\fref$ on the left-hand side, one finds by defining the values $a_k := (f_k - \fref)/\omega$ for all $k \in \N{}$ that
  \bequation\label{eq.3/2_a}
    \underbrace{\frac{f_k - \fref}{\omega}}_{a_k} - \underbrace{\frac{f_{k+1} - \fref}{\omega}}_{a_{k+1}} \geq \underbrace{\frac{(f_{k+1} - \fref)^{3/2}}{\omega^{3/2}}}_{a_{k+1}^{3/2}}.
  \eequation
  One obtains from this inequality that $a_k \geq a_{k+1}^{3/2}$, which when $a_k \geq 1$ (which is to say that $f_k - \fref \geq \omega = \zeta^2/\kappa^3$) gives \eqref{eq.cubic_region11_1}.  Otherwise, \eqref{eq.3/2_a} also yields
  \bequationNN
    \baligned
      \frac{1}{a_{k+1}^{1/2}} - \frac{1}{a_k^{1/2}}
        &\geq \frac{1}{a_{k+1}^{1/2}} - \frac{1}{(a_{k+1} + a_{k+1}^{3/2})^{1/2}} \\
        &= \frac{(a_{k+1} + a_{k+1}^{3/2})^{1/2} - a_{k+1}^{1/2}}{a_{k+1}^{1/2}(a_{k+1} + a_{k+1}^{3/2})^{1/2}} = \frac{(1 + a_{k+1}^{1/2})^{1/2} - 1}{a_{k+1}^{1/2}(1 + a_{k+1}^{1/2})^{1/2}}.
    \ealigned
  \eequationNN
  The right-hand side above is a monotonically decreasing function of $a_{k+1}^{1/2}$ over $a_{k+1} \in (0,1]$.  Hence, when $a_k < 1$ (which is to say that $f_k - \fref < \omega = \zeta^2/\kappa^3$), which implies that $a_{k+1} < 1$, one finds from the above that
  \bequation\label{eq.recall_later}
    \frac{1}{\sqrt{a_{k+1}}} \geq \frac{1}{\sqrt{a_k}} + \frac{\sqrt{2} - 1}{\sqrt{2}}.
  \eequation
  Rearranging this inequality, one obtains \eqref{eq.cubic_region11_2}.
    
  If, with $x_{k+m} \in \region{1}$, an algorithm offers \eqref{eq.RC_region1_alg2_mstep}, then the desired conclusions hold using the same arguments above with \eqref{eq.RC_region1_alg2_mstep} in place of \eqref{eq.RC_region1_alg2}.
  \qed
\eproof

With Theorems~\ref{th.RC_region1_alg1} and \ref{th.RC_region1_alg2}, we can characterize the behavior at points in $\region{1}$ of our algorithms from \S\ref{sec.algorithms}.  This is captured in the following corollary.  (For the results for \RG{} and \RGA{} in this corollary, only Assumption~\ref{ass.first} is needed.  We invoke Assumption~\ref{ass.second} for the sake of being concise as it is needed for the other methods.)

\bcorollary\label{cor.region1}
  Suppose Assumptions~$\ref{ass.first}$ and $\ref{ass.second}$ hold.  Then, the following hold true.
  \benumerate
    \item[$($a$)$] For the \RG{} method, inequality~\eqref{eq.RC_region1_alg1} holds for all $k \in \N{}$ with $\zeta = 2l_1$.
    \item[$($b$)$] For the \RGA, \TRG, and \TRH{} methods, inequality~\eqref{eq.RC_region1_alg1_mstep} holds for all $k \in \N{}$ with $\zeta \in \R{}_{>0}$ and $m \in \N{}$ both sufficiently large relative to functions that depends on $L_1$ and the algorithm parameters but are independent of $k$.
    \item[$($c$)$] For the \RN{} method, inequality~\eqref{eq.RC_region1_alg2} holds for all $k \in \N{}$ with $\zeta \in \R{}_{>0}$ sufficiently large relative to $l_2$.
    \item[$($d$)$] For the \RNA{} method, inequality \eqref{eq.RC_region1_alg2_mstep} holds for all $k \in \N{}$ with $\zeta \in \R{}_{>0}$ and $m \in \N{}$ both sufficiently large relative to functions that depend on $L_2$ and the algorithm parameters but are independent of $k$.
  \eenumerate
  Hence, Theorem~$\ref{th.RC_region1_alg1}$ reveals behavior of the \RG, \RGA, \TRG, and \TRH{} methods, whereas Theorem~$\ref{th.RC_region1_alg2}$ reveals behavior of the \RN{} and \RNA{} methods.
\ecorollary
\bproof
  The fact for \RG{} that \eqref{eq.RC_region1_alg1} holds with $\zeta = 2l_1$ follows from Lipschitz continuity of $g$, the fact that $l_1 > L_1$, and the resulting well-known inequality
  \bequationNN
    f_{k+1} \leq f_k + g_k^Ts_k + \frac{l_1}{2} \|s_k\|^2\ \ \text{for all}\ \ k \in \N{}.
  \eequationNN
  Plugging in $s_k = -g_k/l_1$ and rearranging yields \eqref{eq.RC_region1_alg1}.  As for \RGA{}, for $k=0$ and any $k \in \N{}$ such that $s_{k-1}$ was accepted, one finds that $\nu_k \in [\nu_{\min},\nu_{\max}]$, and, for any $k \in \N{}$ such that~$s_k$ is rejected, one finds $\nu_{k+1} \gets \psi \nu_k$.  These facts, along with the fact that the step will be accepted if $\nu_k > L_1$, implies that \eqref{eq.RC_region1_alg1_mstep} holds for some sufficiently large $\zeta$ and $m$, as claimed.  Similarly, for \TRG{} and \TRH, the desired conclusions can be derived from~\cite{CurtLubbRobi18}; specifically, see Lemmas~2.5 and~2.6 in \cite{CurtLubbRobi18} and note that the trust region radius update implies that an accepted step computed with $\delta_k \equiv \|g_k\|/\nu_k$ leads to $f_k - f_{k+1} \geq \|g_k\|^2/\zeta$ while an accepted step computed with $\delta_k \equiv (\lambda(H_k))_-/\nu_k$ leads to $f_k - f_{k+1} \geq (\lambda(H_k))_-^3/\zeta \geq \|g_k\|^2/\zeta$.
  
  That inequality~\eqref{eq.RC_region1_alg2} holds as stated for the \RN{} method follows as in \cite[Eq.~(4.10)]{NestPoly06}.  That~\eqref{eq.RC_region1_alg2_mstep} holds as stated for \RNA{} follows as described in \cite[\S5.2]{NestPoly06}.
  \qed
\eproof

As we discuss in various specific examples in \S\ref{sec.summary}, the theorems that we have proved in this section allow one to characterize the behavior of the algorithms from \S\ref{sec.algorithms} over much of the search spaces for various (potentially nonconvex) functions of interest.  However, rather than ignore points not included in $\region{1}$, we can capture the behavior of algorithms at additional points by defining additional regions based on higher-order derivatives.  We do this for second-order derivatives next.

\section{Second-Order Regions: Points with Negative Curvature Domination}\label{sec.second-order}

Let us now introduce our notion of a second-order region.  For this definition, recall that the second-order necessary conditions for stationarity with respect to a twice continuously differentiable function $f$ are that $g(x)=0$ and $\lambda(H(x)) \geq 0$.

\bdefinition[Region $\region{2}$]\label{def.2}
  \textit{
  For an objective $f : \R{n} \to \R{}$, scalar $\kappa \in (0,L_2]$, and reference objective value $\fref \in [\finf,\infty)$, let
  \bequation\label{eq.curvature-dominated}
    \region{2} := 
    \{x\in\Lcal\setminus \region{1}: 
    (\lambda(H(x))_-^\tau \geq \kappa(f(x) - \fref) \geq 0\  \text{for some}\  \tau \in [1,3] \}.
  \eequation
  Further, let $\subregion{2}{3}$ be the subset of $\region{2}$ such that the inequality in \eqref{eq.curvature-dominated} holds with $\tau = 3$, let $\subregion{2}{2}$ be the subset of $\region{2}\setminus\subregion{2}{3}$ such that the inequality in \eqref{eq.curvature-dominated} holds with $\tau = 2$, and let $\subregion{2}{1} := \region{2}\setminus(\subregion{2}{2}\cup\subregion{2}{3})$ so $\region{2} = \subregion{2}{1} \cup \subregion{2}{2} \cup \subregion{2}{3}$ and $\subregion{2}{1} \cap \subregion{2}{2} = \subregion{2}{1} \cap \subregion{2}{3} = \subregion{2}{2} \cap \subregion{2}{3} = \emptyset$.
  }
\edefinition

The range for the exponent $\tau$ in this definition can again be justified by considering the pitfalls of values outside of $[1,3]$.  In particular, $\tau \in (0,1)$ disproportionately weighs the negative part of the left-most eigenvalue of the Hessian (as part of a measure of second-order stationarity) at points where it is small in magnitude.  On the other hand, as we shall remark in this section, one can achieve $f(x) - f(x+s) = \Omega(\lambda(H(x))_-^3)$ in certain algorithms.

At any point $x \in \Lcal$ with $f(x) > \fref$, it follows from the definition of $\region{2}$ that one must have $\lambda(H(x)) < 0$ with $\|g(x)\|$ small relative to $\lambda(H(x))_-$, which is to say that the norm of the gradient must be relatively small while the left-most eigenvalue of the Hessian must be negative and relatively large in magnitude.  One can speak of a variety of functions such that $\region{1} \neq \Lcal$, yet $\region{1} \cup \region{2} = \Lcal$, or at least functions for which $\region{2} \neq \emptyset$.  Figure~\ref{fig.illustration} shows segments of domains for two functions wherein one finds elements of $\region{2}$ about first-order stationary points.

\bfigure[ht]
  \centering
  \includegraphics[width=0.45\textwidth]{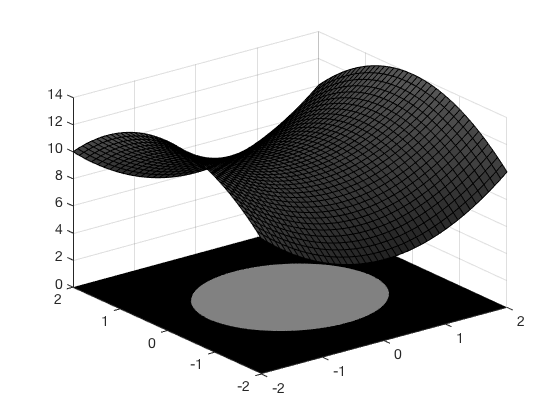}\qquad
  \includegraphics[width=0.45\textwidth]{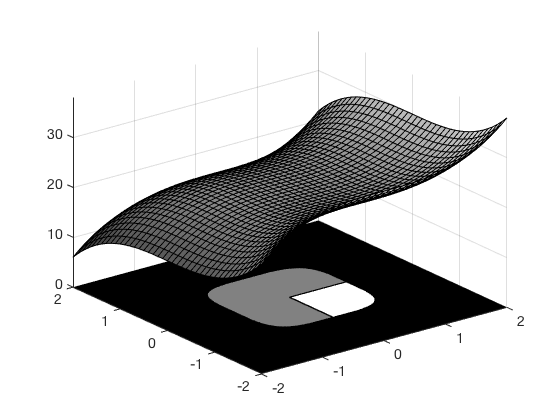}
  \caption{Illustration of $\region{1}$ (black), $\region{2}$ (gray), and $\Lcal\setminus(\region{1}\cup\region{2})$ (white) for the two-dimensional functions $f(x,y) = x^2 - y^2 + 10$ (left) and $f(x) = x^3 - y^3 + 22$ (right).}
  \label{fig.illustration}
\efigure

Given this definition of $\region{2}$, one can provide insight into the performance of an algorithm by tying the reduction obtained with an accepted step to some measure related to the left-most eigenvalue of the Hessian at some point in $\region{2}$.
\smallskip

\setcounter{step}{0}
\noindent
\fcolorbox{black}{gray!20}{\begin{minipage}{0.98\textwidth}
\begin{step}[Region $\region{2}$]\label{step.1_2}
  Attempt to prove that for any accepted step $s_k$ the decrease in the objective function from $x_k$ to $x_{k+1} = x_k + s_k$ satisfies
  \bequation\label{eq.region2_reduction}
    f_k - f_{k+1} = \Omega((\lambda(H(x)))_-^r)\ \ \text{for some}\  x\in\{x_k,x_{k+1}\} \ \text{with} \ x\in\region{2} \ \text{and} \  r > 0.
  \eequation
  If such $(x,r)$ exists, then one can combine \eqref{eq.region2_reduction} and~\eqref{eq.curvature-dominated} to prove a reduction in the objective gap to~$\fref$, i.e., an upper bound for $f_{k+1} - \fref$ as a function of $f_k - \fref$.
  \end{step}
\end{minipage}}
\smallskip

It is implicit in~\eqref{eq.region2_reduction} that one considers the performance of an algorithm over $\region{2}$ only when $\{x_k,x_{k+1}\} \cap \region{2} \neq \emptyset$. This is reasonable since this is precisely when the size of $(\lambda(H(\cdot)))_-$ at $x_k$ and/or $x_{k+1}$ gives information about the size of a potential reduction in the objective through the inequality~\eqref{eq.curvature-dominated} that defines $\region{2}$. 

Naturally, an algorithm should use (approximate) second-order derivative information in order to attain good performance over $\region{2}$.  To demonstrate the instruction above, which we call Step~\ref{step.1_2}--$\region{2}$, we prove the following theorem, which will be useful for characterizing the performance of methods \TRH{}, \RN{}, and \RNA{}.

\btheorem\label{th.RC_region2_alg1}
  Suppose Assumptions~$\ref{ass.first}$ and $\ref{ass.second}$ hold.
  For any algorithm such that having $x_k \in \region{2}$ implies that~\eqref{eq.region2_reduction} holds with $x = x_k$ and $r = 3$ in that
  \bequation\label{eq.RC_region2_alg1}
    f_k - f_{k+1} \geq \frac{1}{\zeta} (\lambda(H_k))_-^3\ \ \text{for some}\ \ \zeta \in [L_2,\infty),
  \eequation
  the following statements hold true.
  \benumerate
    \item[$($a$)$] If $x_k \in \subregion{2}{3}$, then $\{f_k - \fref\}$ decreases as in a linear rate; specifically,
    \bequation\label{eq.RC_region2_linear}
      f_{k+1} - \fref \leq \(1 - \frac{\kappa}{\zeta}\) (f_k - \fref)\ \ \text{where}\ \ \frac{\kappa}{\zeta} \in (0,1].
    \eequation
    \item[$($b$)$] If $x_k \in \subregion{2}{2}$, then it must be true that $\kappa(f_k - \fref) < 1$ and it follows that $\{f_k - \fref\}$ decreases as in a sublinear rate; specifically,
    \bequation\label{eq.RC_region22_sublinear}
      f_{k+1} - \fref \leq \(1 - \(\frac{\kappa^{3/2}}{\zeta}\)\sqrt{f_k - \fref}\)(f_k - \fref).
    \eequation
    \item[$($c$)$] If $x_k \in \subregion{2}{1}$, then it must be true that $\kappa(f_k - \fref) < 1$ and it follows that $\{f_k - \fref\}$ decreases as in a sublinear rate; specifically,
    \bequation\label{eq.RC_region2_sublinear}
      f_{k+1} - \fref \leq \(1 - \frac{\kappa^3}{\zeta}(f_k - \fref)^2\) (f_k - \fref).
    \eequation
  \eenumerate
  Similarly, for any algorithm such that having $x_k \in \region{2}$ implies that
  \bequation\label{eq.RC_region2_alg1_mstep}
    f_k - f_{k+m} \geq \frac{1}{\zeta} (\lambda(H_k))_-^3\ \ \text{for some}\ \ \zeta \in [L_2,\infty)\ \ \text{and}\ \ m \in \N{}
  \eequation
  with $m$ independent of $k$, then $($a$)$, $($b$)$, and $($c$)$ hold with $f_{k+1}$ replaced by $f_{k+m}$.
\etheorem
\bproof
  If $x_k \in \subregion{2}{3}$, then, with \eqref{eq.RC_region2_alg1}, it follows that
  \bequationNN
    f_k - f_{k+1} \geq \frac{1}{\zeta}(\lambda(H_k))_-^3 \geq \frac{\kappa}{\zeta}(f_k - \fref).
  \eequationNN
  Adding and subtracting $\fref$ on the left-hand side and rearranging gives \eqref{eq.RC_region2_linear}.
    
  If $x_k \in \subregion{2}{2}$, which is to say that $(\lambda(H_k))_-^2 \geq \kappa(f_k - \fref)$ while $(\lambda(H_k))_-^3 < \kappa(f_k - \fref)$, then it must be true that $\kappa(f_k - \fref) < 1$.  With \eqref{eq.RC_region2_alg1},
  \bequationNN
    f_k - f_{k+1} \geq \frac{1}{\zeta}(\lambda(H_k))_-^3 \geq \omega^{-1/2} (f_k - \fref)^{3/2}\ \ \text{where}\ \ \omega := \frac{\zeta^2}{\kappa^3}.
  \eequationNN
  Adding and subtracting $\fref$ on the left-hand side, one finds by defining the values $a_k := (f_k - \fref)/\omega = \kappa(f_k - \fref)(\kappa/\zeta)^2 \in [0,1)$ for all $k \in \N{}$ that
  \bequation\label{eq.3232}
    \underbrace{\frac{f_k - \fref}{\omega}}_{a_k} - \underbrace{\frac{f_{k+1} - \fref}{\omega}}_{a_{k+1}} \geq \underbrace{\frac{(f_k - \fref)^{3/2}}{\omega^{3/2}}}_{a_k^{3/2}}.
  \eequation
  One finds from this inequality that $a_{k+1} \leq (1 - \sqrt{a_k})a_k$, which gives \eqref{eq.RC_region22_sublinear}.
    
  If $x_k \in \subregion{2}{1}$, which is to say that $(\lambda(H_k))_- \geq \kappa(f_k - \fref)$ while $(\lambda(H_k))_-^2 < \kappa(f_k - \fref)$, then it must be true that $\kappa(f_k - \fref) < 1$.  In this case, from \eqref{eq.RC_region2_alg1},
  \bequationNN
    f_k - f_{k+1} \geq \frac{1}{\zeta} (\lambda(H_k))_-^3 \geq \frac{1}{\omega^2} (f_k - \fref)^3\ \ \text{where}\ \ \omega := \sqrt{\frac{\zeta}{\kappa^3}}.
  \eequationNN
  Adding and subtracting $\fref$ on the left-hand side, one finds by defining the values $a_k := (f_k - \fref)/\omega = \kappa(f_k - \fref)\sqrt{\kappa/\zeta} \in [0,1)$ for all $k \in \N{}$ that
  \bequationNN
    \underbrace{\frac{f_k - \fref}{\omega}}_{a_k} - \underbrace{\frac{f_{k+1} - \fref}{\omega}}_{a_{k+1}} \geq \underbrace{\frac{(f_k - \fref)^3}{\omega^3}}_{a_k^3}.
  \eequationNN
  One finds from this inequality that $a_{k+1} \leq (1 - a_k^2)a_k$, which gives \eqref{eq.RC_region2_sublinear}.
    
  If, with $x_k \in \region{2}$, an algorithm offers \eqref{eq.RC_region2_alg1_mstep}, then the desired conclusions hold using the same arguments above with \eqref{eq.RC_region2_alg1_mstep} in place of \eqref{eq.RC_region2_alg1}.
  \qed
\eproof

We have the following corollary to Theorem~\ref{th.RC_region2_alg1}.  As previously mentioned, we are only able to state a meaningful result for a few of our second-order algorithms.  After all, one cannot guarantee the performance for the \RG, \RGA, and \TRG{} methods at points in $\region{2}$ since, at any $k \in \N{}$ such that $g_k = 0$ yet $\lambda(H_k) < 0$, these methods would produce zero-norm steps and make no further progress.

\bcorollary\label{cor.region2_alg1}
  Suppose Assumptions~$\ref{ass.first}$ and $\ref{ass.second}$ hold.  Then, the following hold true.
  \benumerate
    \item For the \RN{} method, inequality \eqref{eq.RC_region2_alg1} holds for all $k \in \N{}$ with $\zeta \in \R{}_{>0}$ sufficiently large relative to $l_2$.
    \item For the \TRH{} and \RNA{} methods, inequality \eqref{eq.RC_region2_alg1_mstep} holds for all $k \in \N{}$ with $\zeta \in \R{}_{>0}$ and $m \in \N{}$ both sufficiently large relative to functions that depend on $L_2$ and the algorithm parameters but are independent of $k$.
  \eenumerate
  Hence, Theorem~$\ref{th.RC_region2_alg1}$ reveals behavior of \TRH, \RN, and \RNA{}.
\ecorollary
\bproof
  That \eqref{eq.RC_region2_alg1} holds as stated for \RN{}, and \eqref{eq.RC_region2_alg1_mstep} holds as stated for an \emph{accepted} step for \RNA{} follows from the optimality conditions of the subproblem \eqref{eq.RN}; see, e.g., the proof of \cite[Theorem~5.4]{CartGoulToin11a} or \cite[Equation~(5.28)]{CartGoulToin11b}.  That \eqref{eq.RC_region2_alg1_mstep} holds as stated for an \emph{accepted} step for \TRH{} follows from \cite[Lemma~2.5]{CurtLubbRobi18} and by the definition of $\region{2}$.  Finally, the fact that \eqref{eq.RC_region2_alg1_mstep} holds as stated for arbitrary $k$ for \RNA{} and \TRH{} follows from \cite[\S5.2]{NestPoly06} and \cite[Lemma~2.6]{CurtLubbRobi18}, respectively, which argue that the number of \emph{rejected} steps before the first accepted step or between consecutive accepted steps is uniformly bounded independent of $k$.
  \qed
\eproof

Tables~\ref{tab:region1} and \ref{tab:region2} summarize the RC analysis results that we have presented for the algorithms from \S\ref{sec.algorithms}.  We emphasize that these results have not required referencing any function class.  Rather, they offer insight into performance over the generically defined regions $\region{1}$ and $\region{2}$.  Although interesting in their own right, in the next section we show how the results summarized in these tables may be combined to analyze the performance of the algorithms over collections of regions.

\begin{table}[ht!]
\centering
\caption{Objective decreases over region $\region{1} = \subregion{1}{1}\cup\subregion{1}{2}$ where $\Delta f_k := f_k - \fref$.  Each cell indicates the implied rate and the proved upper bound for $\Delta f_{k+1}/\Delta f_k$.}
  \begin{tabular}{|c|c|c|c|c|} \hline
\rowcolor[gray]{.8}
\multicolumn{2}{|c|}{}              &
\multicolumn{1}{|c|}{\RG/\RGA/\TRG} &
\multicolumn{1}{|c|}{\TRH}          &
\multicolumn{1}{|c|}{\RN/\RNA} \\\hline 
\multirow{2}[4]{*}{$\subregion{1}{1}$}                                                                      &
$\Delta f_k \geq \frac{\zeta^2}{\kappa^3}$                                                                  &
\multirow{2}[4]{*}{$\stackrel{\textnormal{\footnotesize Sublinear}}{1 - \frac{\kappa^2}{\zeta}\Delta f_k}$} &
\multirow{2}[4]{*}{$\stackrel{\textnormal{\footnotesize Sublinear}}{1-\frac{\kappa^2}{\zeta}\Delta f_k}$}   &
$\stackrel{\textnormal{\footnotesize Superlinear}}{\(\frac{\zeta^2}{\kappa^3\Delta f_k}\)^{1/3}}$ \bigstrut \\\cline{2-2}\cline{5-5}
\phantom{}                              &
$\Delta f_k < \frac{\zeta^2}{\kappa^3}$ &
\phantom{}                              &
\phantom{}                              &
$\stackrel{\textnormal{\footnotesize Sublinear}}{\(\frac{1}{1+\frac{\kappa^{3/2}}{\zeta}\(\frac{\sqrt{2}-1}{\sqrt{2}}\)\sqrt{\Delta f_k}  }\)^2}$ \bigstrut \\\cline{1-5}
\multirow{2}[4]{*}{$\subregion{1}{2}$}                                                                      &
$\Delta f_k \geq \frac{\kappa^3}{\zeta^4}$                                                                  &
\multirow{2}[4]{*}{$\stackrel{\textnormal{\footnotesize Linear}}{1 - \frac{\kappa}{\zeta}}$}                &
\multirow{2}[4]{*}{$\stackrel{\textnormal{\footnotesize Linear}}{1-\frac{\kappa}{\zeta}}$}                  &
$\stackrel{\textnormal{\footnotesize Linear}}{\(\frac{\Delta f_0^{1/4}}{\frac{\kappa^{3/4}}{\zeta} + \Delta f_0^{1/4}}\)}$ \bigstrut \\\cline{2-2}\cline{5-5}
\phantom{}                              &
$\Delta f_k < \frac{\kappa^3}{\zeta^4}$ &
\phantom{}                              &
\phantom{}                              &
$\stackrel{\textnormal{\footnotesize Superlinear}}{\(\frac{\zeta^4\Delta f_k}{\kappa^3}\)^{1/3}}$ \bigstrut \\\hline
\end{tabular}
\label{tab:region1}
%
\vspace*{0.5cm}
\centering
\caption{Objective decreases over region $\region{2} = \subregion{2}{1}\cup\subregion{2}{2}\cup\subregion{2}{3}$ where $\Delta f_k := f_k - \fref$.  Each cell indicates the implied rate and the proved upper bound for $\Delta f_{k+1}/\Delta f_k$.}
  \begin{tabular}{|c|c|c|c|} \hline
\rowcolor[gray]{.8}
\multicolumn{1}{|c|}{}              &
\multicolumn{1}{|c|}{\RG/\RGA/\TRG} &
\multicolumn{1}{|c|}{\TRH} &
\multicolumn{1}{|c|}{\RN/\RNA} \\\hline 
$\subregion{2}{1}$ & --- & $\stackrel{\textnormal{\footnotesize Sublinear}}{1 - \frac{\kappa^3}{\zeta}(\Delta f_k)^2}$ & $\stackrel{\textnormal{\footnotesize Sublinear}}{1 - \frac{\kappa^3}{\zeta}(\Delta f_k)^2}$ \bigstrut \\\hline
$\subregion{2}{2}$ & --- & $\stackrel{\textnormal{\footnotesize Sublinear}}{1 - \(\frac{\kappa^{3/2}}{\zeta}\)\sqrt{\Delta f_k}}$ & $\stackrel{\textnormal{\footnotesize Sublinear}}{1 - \(\frac{\kappa^{3/2}}{\zeta}\)\sqrt{\Delta f_k}}$ \bigstrut \\\hline
$\subregion{2}{3}$ & --- & $\stackrel{\textnormal{\footnotesize Linear}}{1 - \frac{\kappa}{\zeta}}$ & $\stackrel{\textnormal{\footnotesize Linear}}{1 - \frac{\kappa}{\zeta}}$ \bigstrut \\\hline
\end{tabular}
\label{tab:region2}
\end{table}

\section{Complete RC Analyses for First- and Second-Order Methods when Minimizing Gradient- and/or Negative Curvature-Dominated Functions}\label{sec.summary}

One way in which RC analysis results may be compared across various algorithms would be to state bounds as in Theorems~\ref{th.RC_region1_alg1}, \ref{th.RC_region1_alg2}, and \ref{th.RC_region2_alg1} (corresponding to algorithms as stated in Corollaries~\ref{cor.region1} and \ref{cor.region2_alg1}).  Indeed, these have all been written in such a way---indicating the reduction in $\{f_k - \fref\}$ for a given iteration---that makes such comparisons straightforward.  That said, equipped with these results, one can also derive complete worst-case performance bounds for algorithms when employed to minimize a function in a class of interest.  This can be done by fitting together results for different regions as if fitting together the pieces of a puzzle.

In this section, we demonstrate a few such complete worst-case performance results for our algorithms from \S\ref{sec.algorithms}.  Our task is to perform the following, which should be understood as the second step introduced on page~\pageref{enum.steps_preliminary}.
\smallskip

\setcounter{step}{1}
\noindent
\fcolorbox{black}{gray!20}{\begin{minipage}{0.98\textwidth}
\begin{step}\label{step.2}
  For an algorithm and different combinations of regions (or subregions), combine results from Step~1 corresponding to these (sub)regions in order to state complete worst-case complexity bounds for the algorithm when it is employed to minimize an objective function for which the search space is completely covered by the combination of regions.  Such bounds hold immediately for functions from classes for which it has been shown that the search space is covered by the combination of regions.
\end{step}
\end{minipage}}
\smallskip

In order to demonstrate Step~2, we provide complete results for our algorithms from~\S\ref{sec.algorithms} when employed to minimize functions from two related classes of (potentially nonconvex) objective functions, defined as follows.

\bdefinition[$(g,H)$-dominated function of degree $(\tau_1,\tau_2)$]\label{def.gH_dominated}
  \textit{
  A twice continuously differentiable function $f$ is $(g,H)$-dominated of degree $(\tau_1,\tau_2) \in [1,2] \times [1,3]$ over $\Lcal$ if for some constant $\kappa \in (0,\min\{L_1,L_2\}]$ it holds that
  \bequation\label{eq.gH_dominated}
    \max\{\|g(x)\|^{\tau_1},(\lambda(H(x)))_-^{\tau_2}\} \geq \kappa(f(x) - \finf)\ \ \text{for all}\ \ x \in \Lcal.
  \eequation
  }
\edefinition

\bdefinition[gradient-dominated function of degree $\tau$]\label{def.g_dominated}
  \textit{
  A continuously differentiable function $f$ is gradient-dominated of degree $\tau \in [1,2]$ over $\Lcal$ if for some constant $\kappa \in (0,L_1]$ it holds that
  \bequation\label{eq.g_dominated}
    \|g(x)\|^{\tau} \geq \kappa(f(x) - \finf)\ \ \text{for all}\ \ x \in \Lcal.
  \eequation
  }
\edefinition

Observe that if $f$ is twice continuously differentiable and gradient-dominated, then it is also $(g,H)$-dominated since \eqref{eq.gH_dominated} holds with $\tau_1=\tau$ and arbitrary $\tau_2$.  On the other hand, not all $(g,H)$-dominated functions are gradient-dominated.  For concreteness, we provide the following examples for these types of functions.

\bexample
  If $f$ has the property that, at all $x \in \Lcal$, either the gradient norm is large in that $\|g(x)\|^2 \geq \kappa(f(x) - \finf)$ or a direction of sufficiently negative curvature exists in that $\lambda(H(x))_-^3 \geq \kappa(f(x) - \finf)$ for some $\kappa \in \R{}_{>0}$, then $f$ is $(g,H)$-dominated of degree $(2,3)$, meaning $\region{1} \cup \region{2} = \subregion{1}{2} \cup \subregion{2}{3} = \Lcal$.  For example, these properties hold for functions satisfying the strict saddle property from \cite[Assumption~A2]{JinGeNetrKakaJord17}, with their constants ``$\theta$'' and ``$\gamma$'' chosen sufficiently small relative to their constant~``$\zeta$'', at least at points that are not approximately globally optimal.  (One could, of course, modify Definitions~\ref{def.gH_dominated} and \ref{def.g_dominated} so as not to require \eqref{eq.gH_dominated} or \eqref{eq.g_dominated} at points with $f_k - \finf \leq \epsilon_f$ if one is interested in the behavior of an algorithm until this $\epsilon_f$-optimality condition is satisfied, as we are in this section.)
\eexample

\bexample
  If $f$ satisfies the Polyak-\L{}ojasiewicz (PL) condition \cite{Poly63} for some constant $\kappa \in (0,L_1]$ at all $x \in \Lcal$, then it is gradient-dominated of degree 2.  For such a function, $\region{1} = \subregion{1}{2} = \Lcal$.  Such functions do not necessarily have unique minimizers.  However, they do have the property that any stationary point is a global minimizer.  The PL condition holds at all $x \in \Lcal$ when $f$ is \emph{strongly convex}, but this is also true for other functions that are not convex.  We refer the reader to \cite{KariNutiSchm16} for a discussion on the relationship between the PL and other types of conditions that have been employed in the context of analyzing optimization methods, such as the \emph{error bounds}, \emph{essential strong convexity}, \emph{weak strong convexity}, \emph{restricted secant inequality}, and \emph{quadratic growth} conditions.  When~$f$ has a Lipschitz continuous gradient, the PL condition is the weakest of these conditions except for the quadratic growth condition, though these are equivalent when $f$ is convex.
\eexample

\bexample
  If $f$ is convex and has a minimizer $x_*$, then $f$ is gradient-dominated of degree 1 with $\kappa = 1/R$ over the Euclidean ball with radius~$R$ centered at~$x_*$ \cite[Example~1]{NestPoly06}.  For such a function, $\region{1}$ includes this ball centered at $x_*$ and $\subregion{1}{1} \neq \emptyset$ if $f$ does not satisfy the PL condition over this domain.
\eexample

We can prove a variety of interesting worst-case performance results (with reference value $\fref = \finf$) for our algorithms from \S\ref{sec.algorithms} when employed to minimize $(g,H)$-dominated or gradient-dominated functions.  The following theorems and corresponding corollaries represent a few examples, in which our main goal is to provide an upper bound on the cardinality of the set of iteration numbers
\bequationNN
  \Kcal_f(\epsilon_f) := \{k \in \N{} : f_k - \finf > \epsilon_f\}.
\eequationNN
For each part of the following results, one might be able to improve the constants involved in the stated convergence rates; however, for ease of comparison, we state results with some common constants.  Throughout this section, let $\epsilon \in (0,\infty)$ be a fixed scalar value that we shall use as an upper bound for the accuracy tolerance~$\epsilon_f$.

Our first two theorems offer complexity bounds for \TRH, \RN, and \RNA{} when they are employed to minimize $(g,H)$-dominated functions of different degrees.

\btheorem\label{th.complete1}
  Suppose that Assumptions~$\ref{ass.first}$ and $\ref{ass.second}$ hold and that \TRH, \RN, or \RNA{} is employed to minimize an objective function $f$ such that $\Lcal = \subregion{1}{2} \cup \subregion{2}{3}$ for some constant $\kappa \in (0,\min\{L_1,L_2\}]$ and $\fref=\finf$.  For $\zeta \in [\max\{L_1,L_2\},\infty)$ satisfying the conditions in Corollaries~$\ref{cor.region1}$ and $\ref{cor.region2_alg1}$ for these methods, let
  \bequation\label{eq.xi}
    \xi := \bcases \(1 - \frac{\kappa}{\zeta}\) & \text{for \TRH} \\ \max\left\{\(1 - \frac{\kappa}{\zeta}\),\(\frac{(f_0 - \finf)^{1/4}}{\frac{\kappa^{3/4}}{\zeta} + (f_0 - \finf)^{1/4}}\)\right\} & \text{for \RN{} and \RNA.} \ecases
  \eequation
  Then, the sequence $\{f_k - \finf\}$ decreases at a linear rate with constant $\xi \in (0,1)$ as defined in \eqref{eq.xi} in the sense that, for some $m \in \N{}$ independent of $k$,
  \bequation\label{eq.RC_region1_firstbound}
    f_{k+m} - \finf \leq \xi(f_k - \finf)\ \ \text{for all}\ \ k \in \N{}.
  \eequation
  Hence, for these methods and any $\epsilon_f \in (0,\epsilon)$,
  \bequation\label{eq.RC_region1_logbound}
    |\Kcal_f(\epsilon_f)| = \Ocal\(\log\(\frac{f_0 - \finf}{\epsilon_f}\)\).
  \eequation
  
  One can go further if $\epsilon_f \in (0,\kappa^3/\zeta^4) \subseteq (0,1/\max\{L_1,L_2\})$ and there exists some iteration number $\khat \in \N{}$ such that $x_k \in \subregion{1}{2}$ for all $k \geq \khat$.  In this case, \TRH{} offers \eqref{eq.RC_region1_firstbound} and consequently \eqref{eq.RC_region1_logbound}, but the convergence rate for \RN{} and \RNA{} improves to superlinear for $k \geq \khat$.  In particular, assuming without loss of generality that $f_{\khat} - \fref < \kappa^3/\zeta^4$ for all $k \geq \khat$, one finds that \RN{} and \RNA{} yield, for the same $m$ as above,
  \bequation\label{eq.RC_region1_firstbound_Newton}
    0 < \frac{4}{3}\log\(\frac{\kappa^3/\zeta^4}{f_k - \finf}\) \leq \log\(\frac{\kappa^3/\zeta^4}{f_{k+m} - \finf}\)\ \ \text{for all}\ \ k \geq \khat,
  \eequation
  in which case it follows for these methods that
  \bequation\label{eq.RC_region1_logbound_Newton}
    |\Kcal_f(\epsilon_f)| = \Ocal\(\log\(\frac{f_0 - \finf}{\kappa^3/\zeta^4}\)\) + \Ocal\(\log\(\log\(\frac{\kappa^3/\zeta^4}{\epsilon_f}\)\)\).
  \eequation
  
  Finally, if for any of these methods $($i.e., \TRH, \RN, or \RNA$)$ a subsequence of the iterate sequence $\{x_k\}$ converges to $x_*$ with $g(x_*) = 0$ and $\lambda(H(x_*)) > 0$, then the entire iterate sequence $\{x_k\}$ eventually converges quadratically to $x_*$.
\etheorem
\bproof
  Since $\Lcal = \subregion{1}{2} \cup \subregion{2}{3}$, it follows from Theorems~\ref{th.RC_region1_alg1}(a), \ref{th.RC_region1_alg2}(a), and \ref{th.RC_region2_alg1}(a) along with Corollaries~\ref{cor.region1}(b) and \ref{cor.region2_alg1}, all with $\fref = \finf$, that for \TRH, \RN, and \RNA{} the inequality \eqref{eq.RC_region1_firstbound} holds for some $m \in \N{}$ for the values of~$\xi$ as stated in \eqref{eq.xi}.  (See also Tables~\ref{tab:region1} and \ref{tab:region2}.)  Applying this fact repeatedly, one finds that
  \bequationNN
    \xi^{k/m}(f_0 - \finf) \geq f_k - \finf\ \ \text{for all}\ \ k \in \{m,2m,3m,\dots\} \subseteq \N{}.
  \eequationNN
  It follows from this inequality that such $k$ satisfy $k \not\in \Kcal_f(\epsilon_f)$ if
  \bequationNN
    \baligned
      \xi^{k/m}(f_0 - \finf) \leq \epsilon_f
      \iff && \frac{f_0 - \finf}{\epsilon_f} &\leq \xi^{-k/m} \\
      \iff && \log\(\frac{f_0 - \finf}{\epsilon_f}\) &\leq -\(\frac{k}{m}\)\log(\xi) \\
      \iff && \frac{m}{-\log(\xi)} \log\(\frac{f_0 - \finf}{\epsilon_f}\) &\leq k,
    \ealigned
  \eequationNN
  from which the bound \eqref{eq.RC_region1_logbound} follows.  In the special case that $\epsilon_f \leq \kappa^3/\zeta^4$ and $x_k \in \subregion{1}{2}$ for all $k \geq \khat$, the first part of the sum in \eqref{eq.RC_region1_logbound_Newton} follows using the same argument as above with $\kappa^3/\zeta^4$ in place of $\epsilon_f$.  Then, for all $k \geq \khat$, the fact that the convergence rate for the \RN{} and \RNA{} methods improves to superlinear follows from Theorem~\ref{th.RC_region1_alg2}(a).  In particular, rearranging \eqref{eq.cubic_region12_2} (with $k+1$ generically replaced by $k+m$ for the same $m \in \N{}$ as above) and taking logs yields~\eqref{eq.RC_region1_firstbound_Newton}.  Then, applying this fact repeatedly, one finds that
  \bequationNN
    \baligned
      \(\frac{4}{3}\)^{(k-\khat)/m}\log\(\frac{\kappa^3/\zeta^4}{f_{\khat} - \finf}\) \leq&\ \log\(\frac{\kappa^3/\zeta^4}{f_k - \finf}\) \\
      &\ \text{for all}\ \ k \in \{\khat+m,\khat+2m,\khat+3m,\dots\} \subseteq \N{}.
    \ealigned
  \eequationNN
  It follows from this inequality that such $k$ satisfy $k \notin \Kcal_f(\epsilon_f)$ if
  \bequationNN
    \baligned
      && \log\(\frac{\kappa^3/\zeta^4}{\epsilon_f}\) &\leq \(\frac{4}{3}\)^{(k-\khat)/m}\log\(\frac{\kappa^3/\zeta^4}{f_{\khat} - \finf}\) \\
      \iff && \(\log\(\frac{\kappa^3/\zeta^4}{f_{\khat} - \finf}\)\)^{-1}\log\(\frac{\kappa^3/\zeta^4}{\epsilon_f}\) &\leq \(\frac{4}{3}\)^{(k-\khat)/m} \\
      \iff && \log\(\(\log\(\frac{\kappa^3/\zeta^4}{f_{\khat} - \finf}\)\)^{-1}\log\(\frac{\kappa^3/\zeta^4}{\epsilon_f}\)\) &\leq \(\frac{k - \khat}{m}\) \log\(\frac43\),
    \ealigned
  \eequationNN
  from which the second term in \eqref{eq.RC_region1_logbound_Newton} follows.  Finally, the fact that the convergence rate for \TRH, \RN, and \RNA{} improves to quadratic if a subsequence of iterates converges to a strong minimizer has been shown in the literature; see \cite[Theorem~3]{NestPoly06}, \cite[Corollary~4.10]{CartGoulToin11a}, and \cite[Theorem~4.1]{FanYuan01}.
  \qed
\eproof

\bcorollary\label{cor.complete1}
  If $f$ is $(g,H)$-dominated of degree $(2,3)$, then $\Lcal = \subregion{1}{2} \cup \subregion{2}{3}$ for some constant $\kappa \in (0,\min\{L_1,L_2\}]$ and $\fref=\finf$.  Hence, when employed to minimize such a function, the behavior of \TRH{}, \RN{}, are \RNA{} is captured by Theorem~$\ref{th.complete1}$.
\ecorollary

One finds from Corollary~\ref{cor.complete1} that when minimizing a $(g,H)$-dominated function of degree $(2,3)$, the behavior of the second-order trust region method \TRH{} is often the same as that of the regularized Newton methods \RN{} and \RNA{}.  The only difference occurs in the special case that the accuracy tolerance is low (i.e., below $\kappa^3/\zeta^4$) and the gradient norms are such that $x_k \in \subregion{1}{2}$ for all large $k$.

Let us now state a result that we shall see requires only that the objective satisfies a weaker form of gradient or negative curvature domination.

\btheorem\label{th.complete2}
  Suppose that Assumptions~$\ref{ass.first}$ and $\ref{ass.second}$ hold and that \TRH, \RN, or \RNA{} is employed to minimize an objective function $f$ such that $\Lcal = \region{1} \cup \region{2}$ for some constant $\kappa \in (0,\min\{L_1,L_2\}]$ and $\fref=\finf$.  For $\zeta \in [\max\{L_1,L_2\},\infty)$ satisfying the conditions in Corollaries~$\ref{cor.region1}$ and $\ref{cor.region2_alg1}$ for these methods, let $\xi \in (0,1)$ be defined as in~\eqref{eq.xi}.  Then, the sequence $\{f_k - \finf\}$ initially decreases at a linear rate with constant~$\xi$ until, for some smallest $\kbar \in \N{}$, one finds that $f_{\kbar} - \finf < \max\{1/\kappa,\epsilon_f\}$.  If $\epsilon_f < 1/\kappa$, then, for $k \geq \kbar$, one of the following cases occurs for some $m \in \N{}$.
    \benumerate
      \item[$($a$)$] If $x_k \in \subregion{1}{2} \cup \subregion{2}{3}$ for all $k \geq \khat$ for some smallest $\khat \geq \kbar$, then, as in Theorem~$\ref{th.complete1}$, the sequence $\{f_k - \finf\}$ decreases linearly $($along the lines of \eqref{eq.RC_region1_firstbound}$)$ and, for sufficiently small $\epsilon_f$, might ultimately decrease superlinearly $($along the lines of \eqref{eq.RC_region1_firstbound_Newton}$)$ for \RN{} and \RNA{}.  Specifically, assuming for simplicity that $\khat = \kbar$, the bound \eqref{eq.RC_region1_logbound} holds for all of these methods and, if $x_k \in \subregion{1}{2}$ for all large $k$ and $\epsilon_f \in (0,\kappa^3/\zeta^4)$, the bound \eqref{eq.RC_region1_logbound_Newton} holds for \RN{} and \RNA.  Moreover, for any of these methods $($i.e., \TRH, \RN, and \RNA$)$, if a subsequence of $\{x_k\}$ converges to $x_*$ with $g(x_*) = 0$ and $\lambda(H(x_*)) > 0$, then the convergence rate of the entire sequence $\{x_k\}$ to $x_*$ is ultimately quadratic.
      \item[$($b$)$] If $x_k \in \subregion{1}{2} \cup (\subregion{2}{2} \cup \subregion{2}{3})$ for all $k \geq \khat$ for some smallest integer $\khat \geq \kbar$, then, in the worst case, the sequence $\{f_k - \finf\}$ eventually decreases sublinearly in that for all $k \in \{\khat,\khat+m,\khat+2m,\khat+3m,\dots\}$ one finds
      \bequation\label{eq.complete11_RN_AN_2}
        f_k - \finf \leq \(\frac{1}{1 + \(\frac{k - \khat}{m}\)\frac{\kappa^{3/2}}{\zeta}\(\frac{\sqrt{2}-1}{\sqrt{2}}\)\sqrt{f_{\khat} - \finf}}\)^2 (f_{\khat} - \finf)
      \eequation
      in which case $($without loss of generality assuming $\khat=\kbar$$)$ it follows that
      \bequation\label{eq.complete11_RN_AN_2_follow}
        |\Kcal_f(\epsilon_f)| = \Ocal\(\log\(\frac{f_0 - \finf}{1/\kappa}\)\) + \Ocal\(\frac{1/\kappa}{\sqrt{\epsilon_f}}\).
      \eequation
      \item[$($c$)$] If $x_k \in (\subregion{1}{1} \cup \subregion{1}{2}) \cup (\subregion{2}{2} \cup \subregion{2}{3})$ for all $k \geq \khat$ for some smallest $\khat \geq \kbar$, then the worst case behavior of \TRH{} is worse than that of \RN{} and \RNA.  In particular, for \TRH, it follows for all $k \in \{\khat,\khat+m,\khat+2m,\khat+3m,\dots\}$ that
    \bequation\label{eq.complete11_TRH_2}
      f_k - \finf \leq \(\frac{1}{1 + \(\frac{k - \khat}{m}\)\frac{\kappa^2}{\zeta} (f_{\khat} - \finf)}\) (f_{\khat} - \finf)
    \eequation
    in which case $($without loss of generality assuming $\khat=\kbar$$)$ it follows that
    \bequation\label{eq.complete11_TRH_2_K}
      |\Kcal_f(\epsilon_f)| = \Ocal\(\log\(\frac{f_0 - \finf}{1/\kappa}\)\) + \Ocal\(\frac{1/\kappa}{\epsilon_f}\).
    \eequation
    On the other hand, for \RN{} or \RNA, it follows for such $k$ that \eqref{eq.complete11_RN_AN_2} holds, in which case $($for simplicity assuming $\khat = \kbar$$)$ it follows that \eqref{eq.complete11_RN_AN_2_follow} holds.
      \item[$($d$)$] If $x_k \in \subregion{2}{1}$ for an infinite number of $k \in \N{}$, then the worst case behavior for all of these methods $($i.e., \TRH, \RN, and \RNA$)$ is the same, i.e., one finds that
    \bequation\label{eq.complete11_TRH_1}
      \baligned
        f_k - \finf \leq&\ \(\sqrt{\frac{1}{1 + \(\frac{k - \kbar}{m}\)\frac{\kappa^3}{\zeta} (f_{\kbar} - \finf)^2}}\) (f_{\kbar} - \finf) \\
        &\ \text{for all large}\ \ k \in \{\kbar,\kbar+m,\kbar+2m,\kbar+3m,\dots\},
      \ealigned
    \eequation
    in which case it follows that
    \bequationNN
      |\Kcal_f(\epsilon_f)| = \Ocal\(\log\(\frac{f_0 - \finf}{1/\kappa}\)\) + \Ocal\(\frac{1/\kappa}{\epsilon_f^2}\).
    \eequationNN
    \eenumerate
\etheorem
\bproof
  Since for $x_k \in \subregion{1}{1} \cup \subregion{2}{1} \cup \subregion{2}{2}$ it must be true that $\kappa(f_k - \finf) < 1$, it follows that while $\kappa(f_k - \finf) \geq 1$ one has that $x_k \in \subregion{1}{2} \cup \subregion{2}{3}$.  For such $k \in \N{}$ with $\kappa(f_k - \finf) \geq 1$, it follows as in the proof of Theorem~\ref{th.complete1} that the sequence $\{f_k - \finf\}$ initially decreases at a linear rate with the constant $\xi$ as given in \eqref{eq.xi}.  Hence, for the remainder of the proof, we may assume that $k \geq \kbar$.
  
  For part (a) with $x_k \in \subregion{1}{2} \cup \subregion{2}{3}$ for all $k \geq \khat$, the conclusions follow using essentially the same arguments as in the proof of Theorem~\ref{th.complete1}.  (One need only also account for iterations $k \in \{\kbar,\dots,\khat-1\}$, but of these there is only a finite number due to the definition of $\khat$.  We ignore these iterations also in the remaining cases of the proof since they do not affect the complexity bounds for small $\epsilon_f$.)
  
  For part (b) with $x_k \in \subregion{1}{2} \cup (\subregion{2}{2} \cup \subregion{2}{3})$ for all $k \geq \khat$, it follows from \eqref{eq.cubic_region12_1}, \eqref{eq.cubic_region12_2}, \eqref{eq.RC_region2_linear}, \eqref{eq.RC_region22_sublinear}, and the fact that $\{f_k - \finf\} \to 0$ that eventually the loosest of these bounds for $f_{k+m} - \finf$ is given by that in \eqref{eq.RC_region22_sublinear}.  Hence, with $\omega := \zeta^2/\kappa^3$ and $a_k := (f_k - \finf)/\omega = \kappa(f_k - \finf)(\kappa/\zeta)^2 \in [0,1)$ for all $k \in \N{}$, it follows as in \eqref{eq.3232} that for sufficiently large $k \geq \khat$ one at least finds
  \bequationNN
    a_k - a_{k+m} \geq a_k^{3/2} \geq a_{k+m}^{3/2}.
  \eequationNN
  Thus, using the same argument as in the proof of Theorem~\ref{th.RC_region1_alg2} that lead from inequality~\eqref{eq.3/2_a} to inequality~\eqref{eq.recall_later}, it follows that
  \bequation\label{eq.recalled}
    \frac{1}{\sqrt{a_{k+m}}} \geq \frac{1}{\sqrt{a_k}} + \frac{\sqrt{2}-1}{\sqrt{2}}.
  \eequation
  Applying this result repeatedly, it follows that
  \bequationNN
    \frac{1}{\sqrt{a_k}} \geq \frac{1}{\sqrt{a_{\khat}}} + \frac{k - \khat}{m}\(\frac{\sqrt{2}-1}{\sqrt{2}}\)\ \ \text{for all}\ \ k \in \{\khat,\khat+m,\khat+2m,\khat+3m,\dots\},
  \eequationNN
  which after rearrangement gives the conclusion in \eqref{eq.complete11_RN_AN_2}.
  
  For part (c) with $x_k \in (\subregion{1}{1} \cup \subregion{1}{2}) \cup (\subregion{2}{2} \cup \subregion{2}{3})$ for all $k \geq \khat$, let us consider \TRH{} separately from \RN{} and \RNA.  For \TRH, it follows from \eqref{eq.RC_region1_sublinear}, \eqref{eq.cubic_region12_1}, \eqref{eq.cubic_region12_2}, \eqref{eq.RC_region2_linear}, \eqref{eq.RC_region22_sublinear}, and the fact that $\{f_k - \finf\} \to 0$ that eventually the loosest of these bounds for $f_{k+m} - \finf$ is given by that in \eqref{eq.RC_region1_sublinear}.  From this, it follows with $a_k := \kappa(f_k - \finf) \in (0,1)$ for all $k \in \N{}$ and $\omega := \kappa/\zeta \in (0,1]$ that one at least finds
  \bequationNN
    a_{k+m} \leq (1 - \omega a_k) a_k \implies \frac{1}{a_{k+m}} \geq \frac{1}{a_k(1 - \omega a_k)} = \frac{1}{a_k} + \frac{\omega}{1 - \omega a_k} \geq \frac{1}{a_k} + \omega,
  \eequationNN
  which, after a repeated use, implies
  \bequationNN
    \frac{1}{a_k} \geq \frac{1}{a_{\khat}} + \(\frac{k - \khat}{m}\)\omega\ \ \text{for all}\ \ k \in \{\khat,\khat+m,\khat+2m,\khat+3m,\dots\}.
  \eequationNN
  Rearranging this inequality leads to the conclusion for \TRH{} in \eqref{eq.complete11_TRH_2}.   Now consider the behavior of \RN{} and \RNA.  First, observe that
  \bequationNN
    f_k - \finf \leq \frac{1}{\kappa} \leq \frac{\zeta^2}{\kappa^3}\ \ \text{for all}\ \ k \geq \kbar.
  \eequationNN
  Hence, it follows from \eqref{eq.cubic_region12_1}, \eqref{eq.cubic_region12_2}, \eqref{eq.cubic_region11_2}, \eqref{eq.RC_region2_linear}, \eqref{eq.RC_region22_sublinear}, and $\{f_k - \finf\} \to 0$ that eventually the loosest of these bounds for $f_{k+m} - \finf$ is given by that in either \eqref{eq.cubic_region11_2} or \eqref{eq.RC_region22_sublinear}, which in either case (as seen above with respect to \eqref{eq.RC_region22_sublinear}) leads to \eqref{eq.recalled}.  Hence, as in the proof for part (b), one is led to the conclusion in \eqref{eq.complete11_RN_AN_2}, from which~\eqref{eq.complete11_RN_AN_2_follow} follows.
  
  For part (d), the worst case behavior of all methods is dictated by \eqref{eq.RC_region2_sublinear}, which with $a_k := \kappa(f_k - \finf) \in (0,1)$ for all $k \in \N{}$ and $\omega := \kappa/\zeta \in (0,1]$ offers
  \bequationNN
    a_{k+m} \leq (1 - \omega a_k^2) a_k.
  \eequationNN
  Hence, one finds that
  \bequationNN
    \frac{1}{a_{k+m}^2} \geq \frac{1}{a_k^2(1 - \omega a_k^2)^2} \geq \frac{1}{a_k^2(1 - \omega a_k^2)} = \frac{1}{a_k^2} + \frac{\omega}{1 - \omega a_k^2} \geq \frac{1}{a_k^2} + \omega.
  \eequationNN
  Applying this result repeatedly, it follows that
  \bequationNN
    \frac{1}{a_k^2} \geq \frac{1}{a_{\kbar}^2} + \(\frac{k - \kbar}{m}\)\omega\ \ \text{for all}\ \ k \in \{\kbar,\kbar+m,\kbar+2m,\kbar+3m,\dots\},
  \eequationNN
  which after rearrangement gives \eqref{eq.complete11_TRH_1}.
  \qed
\eproof

\bcorollary\label{cor.complete2}
  If $f$ is $(g,H)$-dominated of degree $(1,1)$, then $\Lcal = \region{1} \cup \region{2}$ for some constant $\kappa \in (0,\min\{L_1,L_2\}]$ and $\fref=\finf$.  Hence, when employed to minimize such a function, the behavior of \TRH{}, \RN{}, and \RNA{} is captured by Theorem~$\ref{th.complete2}$.
\ecorollary

One finds from Theorem~\ref{th.complete2} that, as in Theorem~\ref{th.complete1}, the behavior of \TRH{} is often the same as that of \RN{} and \RNA{} when minimizing $(g,H)$-dominated functions.  The differences only occur when the accuracy tolerance is small and the algorithm lands on gradient-dominated points of any degree $\tau\in[1,2]$ for large $k$.  Let us also observe that a stronger result than in Theorem~\ref{th.complete2} would be obtained if $f$ were, e.g., assumed to be $(g,H)$-dominated of degree $(1,2)$.  Indeed, in such a situation, one would not need to consider the situation in part (d) of the result.

For our remaining results, we consider gradient-dominated functions of different degrees, about which we are also able to prove results about the first-order methods \RG{} and \RGA{}, as well as the second-order method \TRG{}.  For the following theorems, we are able to borrow from the proofs of Theorems~\ref{th.complete1} and \ref{th.complete2}.

\btheorem\label{th.complete3}
  Suppose that Assumptions~$\ref{ass.first}$ and $\ref{ass.second}$ hold and that any of the algorithms from \S$\ref{sec.algorithms}$ is employed to minimize an objective function $f$ such that $\Lcal = \subregion{1}{2}$ for some constant $\kappa \in (0,\min\{L_1,L_2\}]$ and $\fref=\finf$.  For $\zeta \in [\max\{L_1,L_2\},\infty)$ satisfying the conditions in Corollaries~$\ref{cor.region1}$ and $\ref{cor.region2_alg1}$ for these methods, let
  \bequation\label{eq.xi_more_algs}
    \xi := \bcases \(1 - \frac{\kappa}{\zeta}\) & \text{for \RG, \RGA, \TRG, and \TRH} \\ \max\left\{\(1 - \frac{\kappa}{\zeta}\),\(\frac{(f_0 - \fref)^{1/4}}{\frac{\kappa^{3/4}}{\zeta} + (f_0 - \fref)^{1/4}}\)\right\} & \text{for \RN{} and \RNA.} \ecases
  \eequation
  Then, the sequence $\{f_k - \finf\}$ decreases at a linear rate with constant $\xi \in (0,1)$ as defined in \eqref{eq.xi_more_algs} in the sense that, for some $m \in \N{}$ independent of $k$, the lower bound \eqref{eq.RC_region1_firstbound} holds.  Hence, for any $\epsilon_f \in (0,\epsilon)$, the bound \eqref{eq.RC_region1_logbound} holds.  In addition, if $\epsilon_f \in (0,\kappa^3/\zeta^4) \subseteq (0,1/\max\{L_1,L_2\})$, then the convergence rate for \RN{} and \RNA{} improves to superlinear for large $k$ in the sense that \eqref{eq.RC_region1_firstbound_Newton} holds, leading to \eqref{eq.RC_region1_logbound_Newton}.  Finally, if for \TRG{}, \TRH{}, \RN{}, or \RNA{}, a subsequence of the iterate sequence~$\{x_k\}$ converges to $x_*$ with $g(x_*)=0$ and $\lambda(H(x_*)) > 0$, then the entire sequence~$\{x_k\}$ eventually converges quadratically to~$x_*$.
\etheorem
\bproof
  From Theorems~\ref{th.RC_region1_alg1}, \ref{th.RC_region1_alg2}, and \ref{th.RC_region2_alg1} along with Corollaries~\ref{cor.region1} and \ref{cor.region2_alg1}, all with $\fref = \finf$, the conclusions of the theorem follow using the same arguments as in the proof of Theorem~\ref{th.complete1}.  In addition, the fast local convergence rate for \TRG{} under the stated conditions has been proved as \cite[Theorem~4.1]{FanYuan01}.
\eproof

\bcorollary\label{cor.complete3}
  If $f$ is gradient-dominated of degree 2, then $\Lcal = \subregion{1}{2}$ for some constant $\kappa \in (0,\min\{L_1,L_2\}]$ and $\fref=\finf$.  Hence, when employed to minimize such a function, the behavior of \RG{}, \RGA{}, \TRG{}, \TRH{}, \RN{}, and \RNA{} is captured by Theorem~$\ref{th.complete3}$.
\ecorollary

In Theorem~\ref{th.complete3}, we find a setting in which the behavior of all of the methods from \S\ref{sec.algorithms} behave similarly, except that \RN{} and \RNA{} eventually converge superlinearly if the accuracy tolerance is small.  We also find that each of the second-order methods ultimately converges quadratically if a strong minimizer is approached.

\btheorem\label{th.complete4}
  Suppose that Assumptions~$\ref{ass.first}$ and $\ref{ass.second}$ hold and that any of the algorithms from \S$\ref{sec.algorithms}$ is employed to minimize an objective function $f$ such that $\Lcal = \region{1}$ for some constant $\kappa \in (0,\min\{L_1,L_2\}]$ and $\fref=\finf$.  For $\zeta \in [\max\{L_1,L_2\},\infty)$ satisfying the conditions in Corollaries~$\ref{cor.region1}$ and $\ref{cor.region2_alg1}$ for these methods, let $\xi \in (0,1)$ be defined as in \eqref{eq.xi_more_algs}.  Then, the sequence $\{f_k - \finf\}$ initially decreases at a linear rate with constant $\xi$ until, for some smallest $\kbar \in \N{}$, one finds that $f_{\kbar} - \finf < \max\{1/\kappa,\epsilon_f\}$.  If $\epsilon_f < 1/\kappa$, then, for $k \geq \kbar$, one of the following cases occurs for some $m \in \N{}$.
  \benumerate
    \item[$($a$)$] If $x_k \in \subregion{1}{2}$ for all $k \geq \khat$ for some smallest $\khat \geq \kbar$, then, as in Theorem~$\ref{th.complete3}$, the sequence $\{f_k - \finf\}$ decreases linearly $($along the lines of \eqref{eq.RC_region1_firstbound}$)$ and, for sufficiently small $\epsilon_f$, might ultimately decrease superlinearly $($along the lines of \eqref{eq.RC_region1_firstbound_Newton}$)$ for \RN{} and \RNA{}.  Moreover, for \TRG{}, \TRH{}, \RN{}, and \RNA{}, if a subsequence of $\{x_k\}$ converges to $x_*$ with $g(x_*)=0$ and $\lambda(H_*) > 0$, then, for these methods, the convergence rate of the entire sequence $\{x_k\}$ to $x_*$ is ultimately quadratic.
    \item[$($b$)$] If $x_k \in \subregion{1}{1}$ for an infinite number of $k \in \N{}$, then the worst-case behavior of \RG{}, \RGA{}, \TRG{}, and \TRH{} is the same in that \eqref{eq.complete11_TRH_2} holds, leading to \eqref{eq.complete11_TRH_2_K}.  On the other hand, for \RN{} and \RNA{}, one finds that  \eqref{eq.complete11_RN_AN_2} holds, leading to \eqref{eq.complete11_RN_AN_2_follow}.
  \eenumerate
\etheorem
\bproof
  From Theorems~\ref{th.RC_region1_alg1}, \ref{th.RC_region1_alg2}, and \ref{th.RC_region2_alg1} along with Corollaries~\ref{cor.region1} and \ref{cor.region2_alg1}, all with $\fref = \finf$, the conclusions of the theorem follow using the arguments as in the proofs of Theorems~\ref{th.complete1}, \ref{th.complete2}, and \ref{th.complete3}.
\eproof

\bcorollary\label{cor.complete4}
  If $f$ is gradient-dominated of degree 1, then $\Lcal = \region{1}$ for some constant $\kappa \in (0,\min\{L_1,L_2\}]$ and $\fref=\finf$.  Hence, when employed to minimize such a function, the behavior of \RG{}, \RGA{}, \TRG{}, \TRH{}, \RN{}, and \RNA{} is captured by Theorem~$\ref{th.complete4}$.
\ecorollary

\section{Discussion}\label{sec.discussion}

RC analysis has advantages and disadvantages.  For putting these in perspective, let us first recall known worst-case performance bounds for the algorithms in \S\ref{sec.algorithms}, as they are currently stated in the literature; see \cite{BirgGardMartSantToin17,CartGoulToin11b,CurtLubbRobi18,NestPoly06}.  In particular, suppose Assumptions~\ref{ass.first} and \ref{ass.second} hold and, for any pair of constants $(\epsilon_1,\epsilon_2) \in (0,\epsilon) \times (0,\epsilon)$, let
\bequationNN
  \Kcal_1(\epsilon_1) := \{k \in \N{} : \|g_k\| > \epsilon_1\}\ \ \text{and}\ \ \Kcal_2(\epsilon_2) := \{k \in \N{} : \lambda(H_k) < -\epsilon_2\}.
\eequationNN
Then, one finds that
\bequation\label{eq.eps1}
|\Kcal_1(\epsilon_1)| =
\begin{cases}
\Ocal\(\frac{f_0 - \finf}{\epsilon_1^2}\) & \text{for $\RG{}$, $\RGA{}$, $\TRG{}$, and $\TRH{}$,} \\
\Ocal\(\frac{f_0 - \finf}{\epsilon_1^{3/2}}\) & \text{for $\RN{}$ and $\RNA{}$.}
\end{cases}
\eequation
and that
\bequation\label{eq.eps3}
|\Kcal_2(\epsilon_2)| =
\begin{cases}
\infty & \text{for $\RG{}$, $\RGA{}$, and $\TRG{}$,} \\
\Ocal\(\frac{f_0 - \finf}{\epsilon_2^3}\) & \text{for $\TRH{}$, $\RN{}$, and $\RNA{}$.}
\end{cases}
\eequation
%


While the bounds \eqref{eq.eps1}--\eqref{eq.eps3} hold under relatively loose assumptions, the conclusions are often extremely pessimistic.  Take the bound for $\RG{}$ in \eqref{eq.eps1}, for example.  It is based on the conclusion that with $k \in \Kcal_1(\epsilon_1)$ and an accepted step,
\bequationNN
  f_k - f_{k+1} \geq \frac{1}{2l_1} \|g_k\|^2 \geq \frac{\epsilon_1^2}{2l_1};
\eequationNN
i.e., it only uses the fact that the reduction in the objective attained at such an iterate is at least $\Omega(\epsilon_1^2)$, which is extremely conservative for small~$\epsilon_1$!  On the other hand, for many nonconvex functions, the search space includes many points at which the gradient is significantly larger in norm relative to the objective suboptimality---e.g., points in region~$\region{1}$---from which the attained objective reduction can be much more significant than the (squared) accuracy tolerance.

Another observation is that, with respect to attaining approximate first-order stationarity, \eqref{eq.eps1} offers the same bound for the second-order method \TRH{} as it does for the first-order methods \RG{} and \RGA.  This points to the disappointing conclusions that have been drawn for second-order trust region methods in terms of worst-case performance; see, e.g., \cite{CartGoulToin11c}.  However, for many nonconvex functions, the search space includes many points at which the gradient norm and/or negative curvature is significant---e.g., points in $\region{1} \cup \region{2}$.  For such functions, we have seen that RC analysis offers bounds for the trust region method \TRH{} that are often more similar to those for the regularized Newton methods \RN{} and \RNA{}.

These comments highlight one of the main benefits of RC analysis, namely, that it can offer less pessimistic perspectives on the performance of methods when minimizing certain interesting classes of functions.  However, RC analysis does have some disadvantages.  For one thing, towards attempting to tie the reduction $f_k - f_{k+1}$ to the global error between $f_k$ and some limiting value of the objective attained by an algorithm, we have introduced the reference value $\fref$ that might be considered to be strictly larger than the global minimum $\finf$.  This is useful so that one might be able to use regions to describe the search spaces for functions that one might not be able to minimize to global optimality from all starting points.  Alternatively, if one were only to consider $\fref = \finf$, then, e.g., $\region{1}$ might not include points about local minimizers that are not global minimizers.  All of this being said, it should be clear that the introduction of this reference value puts RC~analysis in no worse of a position than a contemporary worst-case analysis focused on an algorithm attaining (approximate) $p$th-order stationarity.  After all, in the most extreme case for, say, $p=1$, one can consider the reference value $\fref$ to be a placeholder for $\sup_{x\in\R{n}} \{f(x) : \|g(x)\| \leq \epsilon_1\}$ for some $\epsilon_1 \in (0,\infty)$ so that $\region{1}$ at least covers some points at which an algorithm seeking (approximate) first-order stationarity would not yet have terminated.  Put another way: An analysis based on attaining $\|g_k\| \leq \epsilon_1$ also might not offer any guarantees about the number of iterations required to obtain an objective value near the global minimum $\finf$.

Let us now discuss ways in which one can go beyond the results in Theorems~\ref{th.complete1}--\ref{th.complete4}.  In particular, using similar proof strategies, one could prove complete complexity bounds for an algorithm employed to minimize other classes of functions.  For example, if for some class of coercive functions---that are not necessarily $(g,H)$-dominated---one has that $x \in \region{1} \cup \region{2}$ for all $x \in \R{n}$ such that $f(x) \geq \fbar$ for some $\fbar \in [\finf,f_0]$, then one can invoke the results of Theorems~\ref{th.complete1}--\ref{th.complete4} to characterize the behavior of an algorithm at least until $f_k \leq \fbar$ for some $k \in \N{}$.  For all remaining $k \in \N{}$, one can invoke a more conservative bound (e.g., from \eqref{eq.eps1}--\eqref{eq.eps3}) or more refined results depending on the behavior of the algorithm about points in the search space with lower objective values.

One could also obtain different types of results by partitioning regions differently.  For example, if desired for potentially stronger results for a particular class of functions, one could partition $\region{1} = \subregion{1}{2} \cup \subregion{1}{\bar\tau}$ where $\subregion{1}{\bar\tau}$ is the largest subset of $\region{1} \setminus \subregion{1}{2}$ such that the inequality in \eqref{eq.gradient-dominated} holds with $\tau = \bar\tau$.  One then could, e.g., include a separate case along the lines in Theorem~\ref{th.RC_region1_alg1} to derive a certain rate of decrease for $x_k \in \subregion{1}{\bar\tau}$.  The same could be done when partitioning~$\region{2}$ as well.

We also remark that one might consider a \emph{gap} left by RC analysis results to motivate the design of modifications to an algorithm, such as to have the algorithm compute a different type of step or modify some feature of the step computation in order to close the gap.  As an example of the former type of motivation, one can again refer to \cite{CarmHindDuchSidf17} in which a negative curvature direction is computed if/when an accelerated gradient descent method is not behaving as it would when applied to minimize a strongly convex function.  This helps such an algorithm escape neighborhoods about negative-curvature-dominated points.  Another example is the method from \cite{CurtRobi17} that chooses between two types of steps (a first- or a second-order step) depending on which offers a larger predicted reduction in the objective.  As for the second type of motivation, one merely need consider our \TRH{} method, which chooses the trust region radius in each step depending on properties of derivative values.  By doing this, we have seen that \TRH{}---more than the similar method \TRG{}---is able to attain some of the nice features of both the first-order methods \RG{} and \RGA, as well as of the second-order methods \RN{} and \RNA.

\section{Higher-Order Regions, Algorithms, and Analysis}\label{sec.higher-order}

Let us now turn to setting out some fundamental concepts to extend RC analysis to scenarios involving higher-order derivatives.  Let us begin by stating the following assumption, which we shall assume to hold throughout this section.  We employ similar notation as used, e.g., in \cite{BirgGardMartSantToin17}; in particular, the $p$th-order derivative of a function~$f$ at~$x$ is given by the $p$th-order tensor $\nabla^p f(x)$, and the application of this tensor $j \in \N{}$ times to a vector~$s \in \R{n}$ is written as $\nabla^p f(x)[s]^j$.

\bassumption\label{ass.higher}
  The function $f : \R{n} \to \R{}$ is $\pbar$-times continuously differentiable and bounded below by $\finf := \inf_{x\in\R{n}} f(x) \in \R{}$.  In addition, over an open convex set~$\Lcal^+$ containing $\Lcal$ and for each $p \in \{1,\dots,\pbar\}$, the $p$th-order derivative of $f$ is bounded in norm by $M_p \in \R{}_{>0}$ and Lipschitz continuous with Lipschitz constant $L_p \in \R{}_{>0}$ in that
  \bequationNN
    \baligned
      \|\nabla^p f(x)\|_{[p]} &\leq M_p\ \ \text{and} \\
      \|\nabla^p f(x) - \nabla^p f(\xbar)\|_{[p]} &\leq (p-1)! L_p \|x - \xbar\|_2\ \ \text{for all}\ \ (x,\xbar) \in \R{n} \times \R{n},
    \ealigned
  \eequationNN
  where $\|\cdot\|_{[p]}$ denotes the tensor norm recursively induced by $\|\cdot\|$; see \cite[eq.~\textnormal{(2.2)--(2.3)}]{BirgGardMartSantToin17}.
\eassumption

Let us now generalize Definitions~\ref{def.1} and \ref{def.2}.  To do this, let us show that the left-hand side values with largest exponents, namely, $\|g(x)\|^2$ and $(\lambda(H(x)))_-^3$, in Definitions~\ref{def.1} and \ref{def.2} are proportional to the reductions one attains by minimizing a regularized function involving $p$th-order derivatives of the objective at $x \in \Lcal$.  Specifically, for each $p \in \{1,\dots,\pbar\}$, let $v_p(x,\cdot) : \R{n} \to \R{}$ represent the sum of the $p$th-order term of a Taylor series approximation of $f$ centered at $x \in \Lcal$ and a $(p+1)$st-order regularization term, i.e., let $v_p(x,\cdot)$ be defined for all $s \in \R{n}$ by
\bequationNN
  v_p(x,s) = \frac{1}{p!} \nabla^p f(x)[s]^p + \frac{1}{p+1}\|s\|^{p+1}.
\eequationNN
This model is coercive, so it has a minimum norm global minimizer $s_{v_p}(x) \in \R{n}$ with which we can define the reduction function $\Delta v_p : \Lcal \to \R{}$ by
\bequationNN
  \Delta v_p(x) = v_p(x,0) - v_p(x,s_{v_p}(x)) \geq 0.
\eequationNN

We claim that an appropriate generalization of Definitions~\ref{def.1} and \ref{def.2} involves
\bequationNN
  \Delta_p(x) := p(p+1)\Delta v_p(x)\ \ \text{for any}\ \ p \in \{1,\dots,\pbar\}.
\eequationNN
In particular, we now introduce the following definition for $\region{p}$ for all $p \in \{1,\dots,\pbar\}$.

\bdefinition[Region $\region{p}$]\label{def.p}
  \textit{
  For an objective $f : \R{n} \to \R{}$, scalar $\kappa \in (0,L_p]$, and reference objective value $\fref \in [\finf,\infty)$, let
  \bequation\label{eq.pth-dominated}
    \region{p}:=
    \{ x\in\Lcal \setminus \region{p-1}: 
    (\Delta_p(x))^\tau \geq \kappa(f(x) - \fref) \geq 0\  \text{for some} \ \tau \in [1,p+1]\}.
  \eequation
  Further, let $\subregion{p}{p+1}$ be the subset of $\region{p}$ such that the inequality in \eqref{eq.pth-dominated} holds with $\tau = p+1$, and recursively for $q \in \{p,p-1,\dots,1\}$ let $\subregion{p}{q}$ be the subset of $\region{p}\setminus(\subregion{p}{p+1}\cup \subregion{p}{p} \cup \cdots \cup \subregion{p}{q+1})$ such that the inequality in \eqref{eq.pth-dominated} holds with $\tau = q$.
  }
\edefinition

This definition is consistent with Definitions~\ref{def.1} and \ref{def.2}, as the following shows.

\blemma\label{th.Delta_m}
  For $\pbar \geq 2$, it follows that, for any $x \in \Lcal$,
  \bequationNN
    \Delta_1(x) = \|g(x)\|^2\ \ \text{and}\ \ \Delta_2(x) = (\lambda(H(x)))_-^3.
  \eequationNN
\elemma
\bproof
  Let $x \in \Lcal$ be arbitrary.  Since $v_1(x,s) = g(x)^Ts + \thalf \|s\|^2$, one finds that the global minimizer of $v_1(x,\cdot)$ is $s_{v_1}(x) = -g(x)$, meaning that
  \bequationNN
    \baligned
      \Delta v_1(x)
        &= v_1(x,0) - v_1(x,s_{v_1}(x)) \\
        &= -g(x)^Ts_{v_1}(x) - \frac{1}{2}\|s_{v_1}(x)\|^2 = \frac{1}{2} \|g(x)\|^2,
    \ealigned
  \eequationNN
  as desired.  Now consider $v_2(x,s) = \thalf s^TH(x)s + \tfrac13 \|s\|^3$.  If $H(x) \succeq 0$, then the minimum norm global minimizer of $v_2(x,\cdot)$ is $s_{v_2}(x) = 0$.  Otherwise, the global minimum of~$v_2(x,\cdot)$ is achieved at an eigenvector $s_{v_2}(x)$ corresponding to the left-most eigenvalue of $H(x)$, scaled so that it satisfies the first-order condition
  \bequationNN
    (H(x) + \|s_{v_2}(x)\| I)s_{v_2}(x) = 0,
  \eequationNN
  which in particular implies that $\|s_{v_2}(x)\| = -\lambda(H(x))$.  Thus,
  \bequationNN
    \baligned
      \Delta v_2(x)
        &= v_2(x,0) - v_2(x,s_{v_2}(x)) \\
        &= -\half s_{v_2}(x)^TH(x)s_{v_2}(x) - \frac{1}{3} \|s_{v_2}(x)\|^3 \\
        &= -\half \lambda(H(x)) \|s_{v_2}(x)\|^2 - \frac{1}{3} \|s_{v_2}(x)\|^3 \\
        &= \frac{1}{2} |\lambda(H(x))|^3 - \frac{1}{3} |\lambda(H(x))|^3 = \frac{1}{6} |\lambda(H(x))|^3.
    \ealigned
  \eequationNN
  Combining the results of the two cases yields the desired conclusion.
  \qed
\eproof

In order to demonstrate RC analysis results pertaining to $\region{p}$, let us consider a $p$th-order extension of \RG{} and \RN{}.  (The method here can be seen as a special case of the AR$p$ method from \cite{BirgGardMartSantToin17}.)  Let the $p$th-order Taylor series approximation of $f$ at $x \in \Lcal$ be denoted as $t_p(x,\cdot) : \R{n} \to \R{}$, which is given by
\bequationNN
  t_p(x,s) = f(x) + \sum_{j=1}^p \frac{1}{j!} \nabla^j f(x) [s]^j.
\eequationNN
We now define the $\Rp$ method as one that, for all $k \in \N{}$, sets $x_{k+1} \gets x_k + s_{w_p}(x_k)$, where $s_{w_p}(x_k)$ is the minimum-norm global minimizer of a regularized Taylor series approximation function $w_p(x,\cdot) : \R{n} \to \R{}$ defined by
\bequationNN
  w_p(x,s) = t_p(x,s) + \frac{l_p}{p+1}\|s\|^{p+1},\ \ \text{where}\ \ l_p \in \(\frac{(p+1)L_p}{p},\infty\).
\eequationNN

One can draw useful conclusions about the behavior of the \Rp{} method by using the following two example results, which parallel Theorems~\ref{th.RC_region1_alg1}, \ref{th.RC_region1_alg2}, and \ref{th.RC_region2_alg1}.  Our first result can be used to analyze the behavior of \Rp{} over $\region{1}$ using a known decrease property related to its gradient at a point after an accepted step; see \cite{BirgGardMartSantToin17}.

\btheorem\label{th.RC_regionp_alg2}
  Suppose Assumption~$\ref{ass.higher}$ holds.  For any algorithm such that $x_{k+1} \in\region{1}$ implies that~\eqref{eq.region1_reduction} holds with $x = x_{k+1}$ and $r = (p+1)/p$ in that
  \bequation\label{eq.RC_regionp_alg2}
    f_k - f_{k+1} \geq \frac{1}{\zeta} \|g_{k+1}\|^{(p+1)/p}\ \ \text{for some}\ \ \zeta \in (0,\infty),
  \eequation
  the following statements hold true.
  \benumerate
    \item[$($a$)$] If $x_{k+1} \in \subregion{1}{2}$ and $f_k - \fref \geq (\kappa^{p+1}/\zeta^{2p})^{1/(p-1)}$, then $\{f_k - \fref\}$ has decreased as in a linear rate in the sense that the following inequality holds:
    \bequation\label{eq.higher_11}
      f_{k+1} - \fref \leq \(\frac{(f_0 - \fref)^{(p-1)/(2p)}}{\frac{\kappa^{(p+1)/(2p)}}{\zeta} + (f_0 - \fref)^{(p-1)/(2p)}}\)(f_k - \fref).
    \eequation
    On the other hand, if $x_{k+1} \in \subregion{1}{2}$ and $f_k - \fref < (\kappa^{p+1}/\zeta^{2p})^{1/(p-1)}$, then the sequence has decreased as in a superlinear rate in the sense that
    \bequation\label{eq.higher_12}
      f_{k+1} - \fref \leq \(\frac{f_k - \fref}{\(\frac{\kappa^{p+1}}{\zeta^{2p}}\)^{1/(p-1)}}\)^{(p-1)/(p+1)} (f_k - \fref).
    \eequation
    \item[$($b$)$] If $x_{k+1} \in \subregion{1}{1}$, then it must be true that $\kappa(f_{k+1} - \fref) < 1$ and there are two cases: If $f_k - \fref \geq \zeta^p/\kappa^{p+1}$, then $\{f_k - \fref\}$ has decreased superlinearly in that
    \bequation\label{eq.higher_21}
      f_{k+1} - \fref \leq \(\frac{\zeta^p}{\kappa^{p+1}(f_k - \fref)}\)^{1/(p+1)}(f_k - \fref),
    \eequation
    whereas, if $f_k - \fref < \zeta^p/\kappa^{p+1}$, then the sequence has decreased as in a sublinear rate in the sense that the following inequality holds:
    \bequation\label{eq.higher_22}
      f_{k+1} - \fref \leq \(\frac{1}{1 + \frac{\kappa^{(p+1)/p}}{\zeta}\(\frac{2^{1/p}-1}{2^{1/p}}\)(f_k - \fref)^{1/p}}\)^p (f_k - \fref).
    \eequation
  \eenumerate
  Similarly, for an algorithm such that having $x_{k+m} \in \region{1}$ implies that
  \bequation\label{eq.RC_regionp_alg2_mstep}
    f_k - f_{k+m} \geq \frac{1}{\zeta} \|g_{k+m}\|^{(p+1)/p}\ \ \text{for some}\ \ \zeta \in (0,\infty)\ \ \text{and}\ \ m \in \N{},
  \eequation
  with $m$ independent of $k$, then $($a$)$ and $($b$)$ hold with 
  $(x_{k+1},f_{k+1})$ replaced by 
  $(x_{k+m},f_{k+m})$.
\etheorem
\bproof
  If $x_{k+1} \in \subregion{1}{2}$, then, with \eqref{eq.RC_regionp_alg2}, it follows that
  \bequationNN
    \baligned
      f_k - f_{k+1} \geq \frac{1}{\zeta} \|g_{k+1}\|^{(p+1)/p} \geq&\ \omega^{(p-1)/(2p)} (f_{k+1} - \fref)^{(p+1)/(2p)} \\
      \text{where}\ &\ \omega := \(\frac{\kappa^{p+1}}{\zeta^{2p}}\)^{1/(p-1)}.
    \ealigned
  \eequationNN
    Adding and subtracting $\fref$ on the left-hand side, one finds by defining the values $a_k := (f_k - \fref)/\omega$ for all $k \in \N{}$ that
    \bequation\label{eq.p+1/p+2_a}
      \underbrace{\frac{f_k - \fref}{\omega}}_{a_k} - \underbrace{\frac{f_{k+1} - \fref}{\omega}}_{a_{k+1}} \geq \underbrace{\frac{(f_{k+1} - \fref)^{(p+1)/(2p)}}{\omega^{(p+1)/(2p)}}}_{a_{k+1}^{(p+1)/(2p)}}.
    \eequation
    One finds from this inequality that
    \bequationNN
      \frac{a_k}{a_{k+1}} \geq 1 + \frac{1}{a_{k+1}^{(p-1)/(2p)}} \geq 1 + \frac{1}{a_0^{(p-1)/(2p)}} \in (1,\infty),
    \eequationNN
    which gives \eqref{eq.higher_11}.  That said, if $a_k < 1$ (which is to say that $f_k - \fref < \omega$), then one finds from \eqref{eq.p+1/p+2_a} that $a_{k+1} \leq a_k^{2p/(p+1)}$, from which \eqref{eq.higher_12} follows.
    
    If $x_{k+1} \in \subregion{1}{1}$, which is to say that $\|g_{k+1}\| \geq \kappa (f_{k+1} - \fref)$ while $\|g_{k+1}\|^2 < \kappa (f_{k+1} - \fref)$, then it must be true that $\kappa(f_{k+1} - \fref) < 1$.  Hence, with \eqref{eq.RC_regionp_alg2},
    \bequationNN
      f_k - f_{k+1} \geq \frac{1}{\zeta}\|g_{k+1}\|^{(p+1)/p} \geq \omega^{-1/p} (f_{k+1} - \fref)^{(p+1)/p}\ \ \text{where}\ \ \omega := \frac{\zeta^p}{\kappa^{p+1}}.
    \eequationNN
    Adding and subtracting $\fref$ on the left-hand side, one finds by defining the values $a_k := (f_k - \fref)/\omega$ for all $k \in \N{}$ that
    \bequation\label{eq.p+1/p_a}
      \underbrace{\frac{f_k - \fref}{\omega}}_{a_k} - \underbrace{\frac{f_{k+1} - \fref}{\omega}}_{a_{k+1}} \geq \underbrace{\frac{(f_{k+1} - \fref)^{(p+1)/p}}{\omega^{(p+1)/p}}}_{a_{k+1}^{(p+1)/p}}.
    \eequation
    One obtains from this inequality that $a_k \geq a_{k+1}^{(p+1)/p}$, which when $a_k \geq 1$ (which is to say that $f_k - \fref \geq \omega = \zeta^p/\kappa^{p+1}$) gives \eqref{eq.higher_21}.  Otherwise, \eqref{eq.p+1/p_a} also yields
    \bequationNN
      \baligned
        \frac{1}{a_{k+1}^{1/p}} - \frac{1}{a_k^{1/p}}
          &\geq \frac{1}{a_{k+1}^{1/p}} - \frac{1}{\(a_{k+1} + a_{k+1}^{(p+1)/p}\)^{1/p}} \\
          &= \frac{\(a_{k+1} + a_{k+1}^{(p+1)/p}\)^{1/p} - a_{k+1}^{1/p}}{a_{k+1}^{1/p}\(a_{k+1} + a_{k+1}^{(p+1)/p}\)^{1/p}} = \frac{\(1 + a_{k+1}^{1/p}\)^{1/p} - 1}{a_{k+1}^{1/p}\(1 + a_{k+1}^{1/p}\)^{1/p}}.
      \ealigned
    \eequationNN
    The right-hand side above is a monotonically decreasing function of $a_{k+1}^{1/p}$ over $a_{k+1} \in (0,1]$.  Hence, when $a_k < 1$ (which is to say that $f_k - \fref < \omega = \zeta^p/\kappa^{p+1}$), which implies that $a_{k+1} < 1$, one finds from the above that
    \bequationNN
      \frac{1}{a_{k+1}^{1/p}} \geq \frac{1}{a_k^{1/p}} + \frac{2^{1/p} - 1}{2^{1/p}}.
    \eequationNN
    Rearranging this inequality, one obtains \eqref{eq.higher_22}.
    
    If, with $x_{k+m} \in \region{1}$, an algorithm offers \eqref{eq.RC_regionp_alg2_mstep}, then the desired conclusions hold using the same arguments above with \eqref{eq.RC_regionp_alg2_mstep} in place of \eqref{eq.RC_regionp_alg2}.
  \qed
\eproof

Now let us turn to the following result for \Rp{}.  Consistent with our definitions in \S\ref{sec.first-order} and \S\ref{sec.second-order}, one may view this result as an example of following Step~\ref{step.1_2}--$\region{p}$.

\btheorem\label{th.RC_regionp_alg1}
  Suppose Assumptions~\ref{ass.first} and \ref{ass.second} hold.  Then, for any algorithm such that having $x_k \in \region{p}$ implies that the reduction in the objective with an accepted step satisfies
  \bequation\label{eq.RC_regionp_alg1}
    f_k - f_{k+1} \geq \frac{1}{\zeta} (\Delta_p(x_k))^{p+1}\ \ \text{for some}\ \ \zeta \in [L_p,\infty),
  \eequation
  the following statements hold true.
  \benumerate
    \item[$($a$)$] If $x_k \in \subregion{p}{p+1}$, then $\{f_k - \fref\}$ decreases as in a linear rate; specifically,
    \bequation\label{eq.RC_regionp_linear}
      f_{k+1} - \fref \leq \(1 - \frac{\kappa}{\zeta}\) (f_k - \fref)\ \ \text{where}\ \ \frac{\kappa}{\zeta} \in (0,1].
    \eequation
    \item[$($b$)$] If $x_k \in \subregion{p}{q}$ for some $q \in \{1,\dots,p\}$, then it must be true that $\kappa(f_k - \fref) < 1$ and it follows that $\{f_k - \fref\}$ decreases as in a sublinear rate; specifically,
    \bequation\label{eq.RC_regionp_sublinear}
      f_{k+1} - \fref \leq \(1 - \frac{\kappa^{(p+1)/q}}{\zeta}(f_k - \fref)^{(p+1-q)/q}\) (f_k - \fref).
    \eequation
  \eenumerate
  Similarly, for any algorithm such that having $x_k \in \region{p}$ implies that
  \bequation\label{eq.RC_regionp_alg1_mstep}
    f_k - f_{k+m} \geq \frac{1}{\zeta} (\Delta_p(x_k))^{p+1}\ \ \text{for some}\ \ \zeta \in [L_p,\infty)\ \ \text{and}\ \ m \in \N{}
  \eequation
  with $m$ independent of $k$, then $($a$)$ and $($b$)$ hold with $f_{k+1}$ replaced by $f_{k+m}$.
\etheorem
\bproof
  If $x_k \in \subregion{p}{p+1}$, then, with \eqref{eq.RC_regionp_alg1}, it follows that
    \bequationNN
      f_k - f_{k+1} \geq \frac{1}{\zeta}(\Delta_p(x_k))^{p+1} \geq \frac{\kappa}{\zeta}(f_k - \fref).
    \eequationNN
    Adding and subtracting $\fref$ on the left-hand side and rearranging gives \eqref{eq.RC_regionp_linear}.
    
    If $x_k \in \subregion{p}{q}$, which is to say that $(\Delta_p(x_k))^q \geq \kappa(f_k - \fref)$ while $(\Delta_p(x_k))^{q+1} < \kappa(f_k - \fref)$, then it must be true that $\kappa(f_k - \fref) < 1$.  In this case, from \eqref{eq.RC_regionp_alg1},
    \bequationNN
      \baligned
        f_k - f_{k+1} \geq \frac{1}{\zeta} (\Delta_p(x_k))^{p+1} \geq&\ \frac{1}{\omega^{(p+1-q)/q}} (f_k - \fref)^{(p+1)/q} \\
        \text{where}\ &\ \omega := \(\frac{\zeta^q}{\kappa^{p+1}}\)^{1/(p+1-q)}.
      \ealigned
    \eequationNN
    Adding and subtracting $\fref$ on the left-hand side, one finds by defining the values $a_k := (f_k - \fref)/\omega = \kappa(f_k - \fref)(\kappa/\zeta)^{q/(p+1-q)} \in [0,1)$ for all $k \in \N{}$ that
    \bequationNN
      \underbrace{\frac{f_k - \fref}{\omega}}_{a_k} - \underbrace{\frac{f_{k+1} - \fref}{\omega}}_{a_{k+1}} \geq \underbrace{\frac{(f_k - \fref)^{(p+1)/q}}{\omega^{(p+1)/q}}}_{a_k^{(p+1)/q}}.
    \eequationNN
    One finds from this inequality that $a_{k+1} \leq (1 - a_k^{(p+1-q)/q})a_k$, which is \eqref{eq.RC_regionp_sublinear}.
    
    If, with $x_k \in \region{p}$, an algorithm offers \eqref{eq.RC_regionp_alg1_mstep}, then the desired conclusions hold using the same arguments above with \eqref{eq.RC_regionp_alg1_mstep} in place of \eqref{eq.RC_regionp_alg1}.
  \qed
\eproof

Observe that the implied sublinear rate in Theorem~\ref{th.RC_regionp_alg1}(b) improves with larger~$q$.  Indeed, with $q=1$ vs.~$q=p$, one finds reduction factors in \eqref{eq.RC_regionp_sublinear} of
\bequationNN
  1 - \frac{\kappa^{p+1}}{\zeta}(f_k - \fref)^p\ \ \text{vs.}\ \ 1 - \frac{\kappa^{(p+1)/p}}{\zeta}(f_k - \fref)^{1/p}.
\eequationNN
For large $p$, the former can be very close to 1 even for relatively large $f_k - \fref$ (near $1/\kappa$), whereas the latter remains closer to zero due to the exponent on $f_k - \fref$.

Going further, one could explore results that suppose that an algorithm attains $f_k - f_{k+1} = \Omega((\Delta_q(x))^\tau)$ for other $q \in \{1,\dots,p\}$ and some $\tau \geq 1$.  Then, one could combine results from different regions to produce complete RC analysis performance results for different function classes of interest whose search spaces are composed of $\{\region{1},\dots,\region{p}\}$, as was done in \S\ref{sec.summary} for $p\in\{1,2\}$.  Of interest in this context might be a generalization of the \TRH{} method that, to compute $s_k$, minimizes a $p$th-order Taylor series approximation of $f$ at $x_k$ subject to a trust region constraint whose radius is given by $\Delta_j(x_k)^{1/j}$, where $j = \arg\max_{q\in\{1,\dots,p\}}\{\Delta_q(x_k)\}$.

\section{Conclusion}\label{sec.conclusion}

We have proposed a strategy for characterizing the worst-case performance of algorithms for solving nonconvex smooth optimization problems.  The strategy is based on a two-step process: first, one analyzes the behavior of an algorithm over regions defined by generic properties of derivative values, and second, one can combine results from different regions to produce complete worst-case performance results, which in turn can offer results for different function classes of interest.  We have shown how this strategy leads to useful characterizations of a few first- and second-order algorithms, and have demonstrated how to extend the strategy to regions defined by, and for algorithms that make use of, higher-order derivatives.

We imagine that our approach for analyzing worst-case complexity can be generalized or adapted to other settings.  The following are some possibilities. $(i)$~While Assumptions~\ref{ass.first}--\ref{ass.higher} require the $p$th-order derivatives of $f$ to be Lipschitz continuous over $\Lcal^+$ for all $p \in \{1,\dots,\pbar\}$ for some $\pbar \in \N{}$, one might consider a more general setting when these derivatives are only H\"older continuous with exponent $\alpha$ not necessarily equal to one; see, e.g., \cite{CartGoulToin17b}. $(ii)$ One might consider nonmonotone methods and settings in which $f$ is extended-real-valued as long as an algorithm can guarantee that, after some number of iterations, a sufficient reduction in the objective is produced.  Indeed, with the flexibility introduced by $m \in \N{}$, this was all that was required for our results. $(iii)$ One might extend our strategy to offer probabilistic results or to analyze stochastic algorithms.  For example, while one is not able to supply a \emph{deterministic} upper bound for \RG{} over $\region{2}$, one can establish \emph{probabilistic} upper bounds by introducing randomization into the starting point or the step computation; see \cite{JinGeNetrKakaJord17,LeeSimcJordRech16}.  As another example, if one is able to ensure that over some number of iterations an algorithm will offer a sufficiently large \emph{expected} reduction in the objective, then generalized forms of our results might involve $f_k - \E_k[f_{k+m}]$ where $\E_k$ denotes the conditional expectation given that the algorithm has reached $x_k$.  Finally, it is conceivable that one can build results based on inequalities such as \eqref{eq.RC_region1_alg1} that are only guaranteed to hold with certain probability \cite{CartSche17}, although we admit this might be a nontrivial extension of our proposed ideas. $(iv)$ An extension of our strategy to nonsmooth $f$ might be based on replacing the measure $\|g(x)\|$ in \eqref{eq.gradient-dominated} in Definition~\ref{def.1} with the norm of a \emph{proximal} step computed at $x \in \Lcal$.  Similarly, one might extend our strategy to constrained optimization if $\|g(x)\|$ is replaced by the norm of a \emph{projected} gradient step.

\bibliographystyle{plain}
\bibliography{nonconvex_complexity}

\end{document}